\documentclass[hidelinks,onefignum,onetabnum]{siamart220329}

\usepackage{amssymb}
\usepackage{amsfonts}
\usepackage{dsfont}
\usepackage{mathrsfs}
\usepackage{graphicx}

\usepackage{enumerate}
\usepackage{hyperref}
\usepackage{indentfirst}
\usepackage{color}
\usepackage{multicol,colordvi,epsfig,float}
\usepackage{pgfplots}
\pgfplotsset{width=9cm,compat=1.15}

\newcommand{\norm}[2]{\left |  {#1}  \right| _{#2}}
\newcommand{\norma}[2]{\left \|  {#1} \right \|_{#2}}
\newcommand{\prodl}[2]{\left (\ {#1} \ , \ {#2} \ \right )}

\newcommand{\D}[1]{\displaystyle{#1}}

\usepackage{verbatim}

\ifpdf
  \DeclareGraphicsExtensions{.eps,.pdf,.png,.jpg}
\else
  \DeclareGraphicsExtensions{.eps}
\fi

\newsiamremark{remark}{Remark}
\newsiamremark{hypothesis}{Hypothesis}
\crefname{hypothesis}{Hypothesis}{Hypotheses}
\newsiamthm{claim}{Claim}

\headers{Convergence of a time discrete scheme for a chemotaxis-consumption model}{F. Guillén-González, A. L. Corrêa Vianna Filho}

\title{Convergence of a time discrete scheme for a chemotaxis-consumption model\thanks{Submitted to the editors \today.
\funding{This work was partially funded by Grant PGC2018-098308-B-I00 (MCI/AEI/FEDER, UE). FGG has also been financed in part by the Grant US-1381261 (US/JUNTA/FEDER, UE) and Grant P20-01120 (PAIDI/JUNTA/FEDER,UE).}}   }

\author{Francisco Guillén-González\thanks{Departament of Partial Differential Equations and Numerical Analysis, Universidad de Sevilla, Seville, Spain. (\email{guillen@us.es}, \email{acorreaviannafilho@us.es}) .}
\and André Luiz Corrêa Vianna Filho\footnotemark[2]}

\usepackage{amsopn}

\begin{document}

\maketitle

\begin{abstract}
  In the present work we propose and study a time discrete scheme for the following chemotaxis-consumption model (for any $s\ge 1$), $$ \partial_t u - \Delta u  = -  \nabla \cdot (u \nabla v), \quad \partial_t v - \Delta v  = - u^s v \quad \hbox{in $(0,T)\times \Omega$,}$$ endowed with isolated boundary conditions and initial conditions, where $(u,v)$ model cell density and chemical signal concentration. The proposed scheme is defined via a reformulation of the model, using the auxiliary variable $z = \sqrt{v + \alpha^2}$ combined with a Backward Euler scheme for the $(u,z)$-problem and a upper truncation of $u$ in the nonlinear chemotaxis and consumption terms. Then, two different ways of retrieving an approximation for the function $v$ are provided. We prove the existence of solution to the time discrete scheme and establish uniform in time \emph{a priori} estimates, yielding the convergence of the scheme towards a weak solution $(u,v)$ of the chemotaxis-consumption model.
\end{abstract}



\section{Introduction}
  
   Chemotaxis is the movement of cells in response to the concentration gradient of a chemical signal. This is a phenomenon that plays an important role in many biological events of practical interest. One of the first models for chemotaxis is attributed to two works of Keller and Segel in 1970 and 1971 \cite{keller1970initiation,keller1971model}, which are also regarded as a development of the work of Patlak \cite{patlak1953random}. Since then, many researchers have been studying this and other related models, which were derived from the Keller-Segel model with the objective of describing different situations involving chemotaxis as well as to obtain models that are more biologically precise. This effort lead to a series of important theoretical developments on chemotaxis models \cite{bellomo2015toward}.
   
   In the present work we focus on a chemotaxis-consumption model. Let $\Omega$ be a bounded domain of $\mathbb{R}^3$
   and let $\Gamma$ be its boundary. The assumptions on $\Omega$ depend on definitions that have not been introduced yet and therefore will be clarified later, in Remark \ref{remark_hypotheses_domain}. Let $u=u(t,x)$ and $v=v(t,x)$ be the density of cell population and the concentration of chemical substance, respectively, defined on $x \in \Omega$ and $t > 0$. This model is governed by the initial-boundary PDE problem
   \begin{equation} \vspace{-1pt}   \label{problema_P}
     \left\{\begin{array}{l}
        \partial_t u - \Delta u  = -  \nabla \cdot (u \nabla v), \quad
        \partial_t v - \Delta v  = - u^s v, \\
        \partial_\eta u \Big |_{\Gamma}  =  \partial_\eta v \Big |_{\Gamma} = 0, \quad
        u(0)  = u^0, \quad v(0) = v^0,
      \end{array}\right.
    \vspace{-1pt} \end{equation} 
    where $\nabla \cdot (u \nabla v)$ is the chemotaxis term and $u^s$ is the consumption rate, with $s \geq 1$. $\partial_\eta u$ denotes the normal derivative of $u$ on the boundary. The functions $u^0$ and $v^0$ are initial conditions satisfying $u^0 \geq 0$ and $v^0 \geq 0$ a.e. in $\Omega$. The biologically meaningful solutions of \eqref{problema_P} are expected to satisfy
    \begin{equation*} \vspace{-1pt}
      \int_{\Omega}{u(t,x) \ dx} = \int_{\Omega}{u^0(x) \ dx}, \ a.e \ t \in (0, \infty),
    \vspace{-1pt} \end{equation*}  
    $u(t,x) \geq 0$ and $\norma{v^0}{L^{\infty}(\Omega)} \geq v(t,x) \geq 0, \ a.e. \ (t,x) \in (0,\infty) \times \Omega$.
    
    Next we introduce some theoretical studies about problem \eqref{problema_P}, beginning with the case $s = 1$. In \cite{tao2012eventual}, for convex and smooth domains, it is established the existence and uniqueness of a regular solution on $2$D domains and the existence of global weak solutions which become regular after a sufficiently large period of time on $3$D domains. Existence and uniqueness of a regular solution is obtained in \cite{tao2011boundedness} on smooth domains of $\mathbb{R}^d$, for $d \geq 1$, under a smallness assumption on the initial signal concentration $v^0$. In \cite{winkler2012global}, a coupled fluid-chemotaxis-consumption model is studied in smooth convex $2$D and $3$D domains, yielding results analogous to those in \cite{tao2012eventual}. Finally, in \cite{jiang2015global} the authors work to extend the results of \cite{winkler2012global} to smooth nonconvex domains.
    
    Problem \eqref{problema_P}, with $s \geq 1$, is studied in \cite{ViannaGuillen2023uniform}. The authors prove the existence and uniqueness of a strong solution in $2$D domains and the existence of global weak solutions in $3$D under minimal assumptions on the boundary of the domain. These weak solutions are obtained through the limit of solutions of adequate truncated models. Additionally, it is shown that in $2$D the solution of \eqref{problema_P} is uniformly bounded in the $L^{\infty}$-norm and, therefore, the truncated model coincides with \eqref{problema_P} if the truncation parameter is greater than $\norma{u}{L^{\infty}(0,\infty;L^{\infty}(\Omega))}$.
    
    As it can be observed in \cite{ViannaGuillen2023uniform}, \cite{tao2012eventual} and \cite{winkler2012global}, the general results on the existence of solution for \eqref{problema_P} rely on the existence of a bounded energy (not necessarily decreasing) based on the cancellation between the chemotaxis and the consumption effects.  
    
    Regarding the numerical approximation of chemotaxis models, although it is a relevant and growing research topic, we still can find chemotaxis models for which there is a relatively low number of numerical studies. This is the case of the chemotaxis-consumption model \eqref{problema_P}. To the best of our knowledge, we can cite two studies, \cite{duarte2021numerical} and \cite{FranciscoGiordano}, about the numerical approximation of \eqref{problema_P}, both just for the case $s = 1$. In \cite{duarte2021numerical} a chemotaxis-Navier-Stokes system is approached via Finite Elements (FE). In fact, by assuming the existence of a sufficiently regular solution, if the initial data of the scheme are small perturbations of the initial data of this regular solution, then optimal error estimates are deduced. The drawback of this result is that the existence of such a regular solution is not clear in general, especially when we consider polyhedral domains, which are broadly used in numerical simulations.
    
    In \cite{FranciscoGiordano}, motivated by the treatment given to the chemorepulsion model with linear production in \cite{guillen2019unconditionally}, several FE schemes are designed to approximate \eqref{problema_P}, with $s = 1$. The authors focus on FE schemes satisfying properties such as conservation of cells, discrete energy law and approximate positivity rather than convergence. In particular, they present a scheme satisfying a discrete energy law that, in 1D domains, yields decreasing energy. Numerical simulations are carried out to compare the performance of the different schemes.
    
    One of the main difficulties of addressing issues concerning the convergence towards weak solutions of \eqref{problema_P} is probably the lack of energy \emph{a priori} estimates for the solutions of the schemes. Even if we consider only time discretizations of \eqref{problema_P}, the task of designing a convergent scheme does not become much easier. This could be attributed to the complex technique needed in order to cancel the chemoattraction and consumption effects. As far as we know, excepting the case of $1$D domains \cite{FranciscoGiordano}, there is no time discrete scheme for \eqref{problema_P} possessing an energy inequality from which one can obtain estimates for the discrete solutions, yielding convergence.
    
    This is why the present work is devoted to propose a time discrete scheme convergent to \eqref{problema_P}. This convergence will be valid in $3$D domains and based on energy estimates. In addition, the scheme will preserve the properties of positivity and conservation of cells. There is evidence that the preservation of the positivity could possibly enhance the performance of the numerical schemes, avoiding spurious oscillations \cite{guillen2022comparison}.
    
    The design of the time discrete scheme is based on the analysis that was carried out in \cite{ViannaGuillen2023uniform}. In \cite{ViannaGuillen2023uniform}, it was convenient to rewrite \eqref{problema_P} in terms of the variable $\sqrt{v+\alpha^2}$, because the test functions involved in obtaining a discrete energy law become simpler. Hence, in the present work we consider the following reformulation of \eqref{problema_P}
    \begin{equation} \vspace{-1pt}  
      \left\{\begin{array}{l}
        \partial_t u - \Delta u  = -  \nabla \cdot (u \nabla (z)^2), \quad
        \partial_t z - \dfrac{\norm{\nabla z}{}^2}{z} - \Delta z  = - \dfrac{1}{2} u^s 
        \left(z - \dfrac{\alpha^2}{z} \right) , \\
        \partial_\eta u \Big |_{\Gamma}  =  \partial_\eta z \Big |_{\Gamma} = 0, \quad
        u(0)  = u^0, \quad z(0) = \sqrt{v^0 + \alpha^2},
      \end{array}\right.
      \label{problema_P_u_z} 
    \vspace{-1pt} \end{equation} 
    where $\alpha > 0$ is a fixed real number to be chosen later in Lemmas \ref{lemma_estimativa_u_v_m_1_s_intermediario} and \ref{lemma_estimativa_u_v_m_1_s_geq_2}. Since it is proved in \cite{ViannaGuillen2023uniform} that the $v$-equation of \eqref{problema_P} is satisfied in the strong sense, with $v \in L^2(0,\infty;H^2(\Omega))$, one can check by straightforward calculations that \eqref{problema_P_u_z} is equivalent to \eqref{problema_P} if we use the change of variables $z = \sqrt{v + \alpha^2}$. We summarize this statement in the following lemma for further use.
    
    \begin{lemma}
      Problems \eqref{problema_P} and \eqref{problema_P_u_z} are equivalent. More precisely, $(u,z)$ is a weak solution of \eqref{problema_P_u_z} if, and only if, $(u,v)$ is a weak solution of \eqref{problema_P}, with $v = z^2 - \alpha^2$.
      \label{lema_equivalencia_problemas}
    \end{lemma}
    
    For the time discretization, we will divide the interval $[0,\infty)$ in subintervals denoted by $I_n = (t_{n-1}, t_n)$, with $t_0 = 0$ and $t_n = t_{n-1} + k$, where $k > 0$ is the length of the intervals $I_n$. If $\{ z^n \}_n$ is a sequence of functions, then we use the notation
    \begin{equation} \vspace{-1pt}  
      \delta_t z^n = \frac{z^n - z^{n-1}}{k}, \quad \forall n \geq 1,
    \vspace{-1pt} \end{equation} 
    for the discrete time derivative. We will also use the following upper truncation of $u$: for each fixed $m > 0$ we define the function $T^m: \mathbb{R} \rightarrow \mathbb{R}$ given by
    \begin{equation} \vspace{-1pt}  
      T^m(u) = \left \{
      \begin{array}{cl}
        u, & \mbox{ if } u \leq m, \\
        m, & \mbox{ if } u \geq m.
      \end{array}
      \right .
      \label{truncamento_limitado_da_identidade}
    \vspace{-1pt} \end{equation} 
    
    In this paper, we propose the following time discrete scheme:
    
    \noindent  {\bf Initialization: } Consider the initial conditions $u^0 \in L^2(\Omega)$, $z^0 = \sqrt{v^0 + \alpha^2} \in L^{\infty}(\Omega)$ and $v^0 \in L^{\infty}(\Omega)$ introduced above.
    
    \noindent  {\bf Step $n$: }    
    Given non-negative functions $u^{n-1} \in L^2(\Omega)$, $z^{n-1} \in L^{\infty}(\Omega)$ and $v^{n-1} \in L^{\infty}(\Omega)$, 
    \begin{enumerate}
      \item Find  $(u^n,z^n) \in H^2(\Omega)^2$, satisfying the bounds 
      $$u^n(x) \geq 0\quad \hbox{and}\quad \norma{z^{n-1}}{L^{\infty}(\Omega)} \geq z^n(x) \geq \alpha \quad \hbox{a.e. $x \in \Omega$,}
      $$
      and the boundary-value problem 
      \begin{equation} \vspace{-1pt}  \label{NLTD}
        \left\{
        \begin{array}{l}
          \delta_t u^n - \Delta u^n =  \nabla \cdot \Big (  T^m(u^n) \nabla (z^n)^2 \Big ), \\
          \delta_t z^n - \dfrac{\norm{\nabla z^n}{}^2}{z^n} - \Delta z^n = - \dfrac{1}{2} T^m(u^n)^s \left(z^n - \dfrac{\alpha^2}{z^n} \right), \\
          \partial_{\eta} u^n \Big |_{\partial \Omega} = \partial_{\eta} z^n \Big |_{\partial \Omega} = 0.
        \end{array}
        \right.
      \vspace{-1pt} \end{equation} 
      \item Two variants for the approximation of $v$ are possible (equally denoted), either depending on $z^n$ or $u^n$:
      \begin{itemize}
        \item  Find $v^n = v^n(z^n) \in H^2(\Omega)$ as
        \begin{equation} \vspace{-1pt}  \label{v-z}
          v^n = (z^n)^2 - \alpha^2 .
        \vspace{-1pt} \end{equation} 
        \item Find $v^n = v^n(u^n) \in H^2(\Omega)$ as the unique solution of the linear problem
        \begin{equation} \vspace{-1pt}  \label{v-u}
          \delta_t v^n - \Delta v^n + T^m(u^n)^s v^n = 0,
          \quad 
          \partial_{\eta} v^n \Big |_{\partial \Omega} =0.
        \vspace{-1pt} \end{equation} 
      \end{itemize}
    \end{enumerate}

\

    \begin{remark} \label{remark_dependencia_esquema_em_m}
      Note that $u^n$, $z^n$ and $v^n$ depend on $m$. For simplicity, from now on, we consider that $m \in \mathbb{N}$.
    \end{remark}

\

    Using the functions $u^n$, $z^n$ and $v^n$ introduced above in step $n$, we define the piecewise function $u_m^k$ and the locally linear and globally continuous function $\tilde{u}_m^k$ by
    \begin{equation} \vspace{-1pt}  \label{funcao_u^kr_u^k}
      \begin{array}{c}
        u_m^k(t,x) = u^n(x) \mbox{ and } \\[6pt]
        \tilde{u}_m^k(t,x) = u^n(x) + \dfrac{(t - t_n)}{k} \big ( u^n(x) - u^{n-1}(x) \big ), \mbox{ if } t \in [t_{n-1},t_n).
      \end{array}
    \vspace{-1pt} \end{equation} 
    Analogously, we define the functions $z_m^k$, $\tilde{z}_m^k$, $v_m^k$ and $\tilde{v}_m^k$.
    
    Now, we are in position to present the main result that will be proved along the present work.
    \begin{theorem} \label{teo_principal}
      For each $n \in \mathbb{N}$, there exists at least one solution $(u^n,z^n)$  of \eqref{NLTD}, that jointly to $v^n$ defined by \eqref{v-z} or \eqref{v-u} leads us to $(u^n,v^n)$ satisfying
      $$u^n(x) \geq 0, \qquad \norma{v^0}{L^{\infty}(\Omega)} \geq v^n(x) \geq 0 \quad a.e. \ x \in \Omega.$$
      Moreover, up to a subsequence, $(u_m^k,v_m^k)$ converges towards a weak solution $(u,v)$ of \eqref{problema_P} as $(m,k) \to (\infty,0)$.
    \end{theorem}

\

    \begin{remark}
      The number $\alpha$ is a sufficiently small positive real number that is chosen in Lemmas \ref{lemma_estimativa_u_v_m_1_s_intermediario} and \ref{lemma_estimativa_u_v_m_1_s_geq_2} independently of $m$ and $k$. The convergence result given in Theorem \ref{teo_principal} as $(m,k) \to (\infty,0)$ is unconditional, that is, there is not any constraint over $m$ and $k$ as long as $m \to \infty$ and $k \to 0$.
    \end{remark}

\
    
    In particular,  we also prove the  result on existence of weak solutions to \eqref{problema_P} in $3$D domains given in \cite[Theorem 1]{ViannaGuillen2023uniform} but, this time, as a consequence of the convergence of the time discrete scheme introduced in this paper.
    
    \
    
    In $2$D domains, there exists a unique strong solution of \eqref{problema_P} (see \cite{ViannaGuillen2023uniform}). The proof is achieved through the obtaining of stronger $m$-independent estimates for the solution of an adequate truncated problem. Unfortunately, it is not clear how we could adapt these strong estimates for the time discrete scheme. Consequently, in $2$D, the convergence of the whole sequence of solutions of the time discrete scheme \eqref{NLTD} towards the unique strong solution of \eqref{problema_P} as $(m,k) \to (0,\infty)$ remains as an open problem.
    
    \
    
    The rest of the paper is organized as follows. In Section \ref{sec:preliminary} we give some technical results that will be used throughout the paper. In Section \ref{Sec:u-z} we establish the existence of solution $(u^n,z^n)$ of the $(u,z)$-scheme \eqref{NLTD}, some pointwise estimates independent of $(m,k,n)$ and an energy inequality for $(u^n,z^n)$. In Section \ref{subsec:passage_limit_m_to_infty}, starting from this energy inequality, we deduce additional \emph{a priori} estimates for $(u^n,z^n)$, independent of $(m,k,n)$, that allow us to pass to the limit as $(m,k) \to (\infty,0)$, obtaining convergence of \eqref{NLTD} towards the $(u,z)$-problem \eqref{problema_P_u_z}. Finally, in Section \ref{subsec:convergencia_de_v_m^kr} we prove the convergence of $(u^n,v^n)$, with $v^n$ defined by \eqref{v-z} or \eqref{v-u}, towards the $(u,v)$-problem \eqref{problema_P}.

\



\section{Preliminary Results}
  
  \label{sec:preliminary}
  
  For all $g \in L^1(\Omega)$, define
  \begin{equation*} \vspace{-1pt}
    g^{\ast} = \dfrac{1}{\norm{\Omega}{}} \D{\int_{\Omega} g(x) \ dx}
  \vspace{-1pt} \end{equation*}  
  \begin{lemma}[\bf Poincare's Inequality, \cite{Evans2010}]
    There is a constant $C > 0$  such that
    \begin{equation*} \vspace{-1pt}
      \norma{v - v^{\ast}}{W^{1,p}(\Omega)} \leq C \norma{\nabla v}{L^p(\Omega)},
      \quad 
      \forall\, v \in W^{1,p}(\Omega).
    \vspace{-1pt} \end{equation*}  
    \label{lema_desigualdade_poincare_media_nula}
  \end{lemma}
  
  Concerning the regularity of solutions of the Poisson-Neumann problem
  \begin{equation} \vspace{-1pt}  
    \left \{ \begin{array}{rl}
      - \Delta z + z & = f  \quad \mbox{ in } \Omega, \\
      \partial_\eta z \Big |_{\Gamma} & = 0,
    \end{array} \right.
    \label{Neumann_problem}
  \vspace{-1pt} \end{equation} 
  we present the following.
  
  \begin{definition}
    Let $z \in H^1(\Omega)$ be a weak solution of \eqref{Neumann_problem} with $f \in L^p(\Omega)$. If this implies that $z \in W^{2,p}(\Omega)$ with
    \begin{equation*} \vspace{-1pt}
      \norma{z}{W^{2,p}(\Omega)} \leq C \norma{- \Delta z + z}{L^p(\Omega)},
    \vspace{-1pt} \end{equation*}  
    then we say that the Poisson-Neumann problem \eqref{Neumann_problem} has $W^{2,p}$-regularity. In the hilbertian case $p = 2$ we say $H^2$-regularity.
    \label{defi_regularidade_H_m}
  \end{definition}

\

  \begin{lemma}
    Let $\Omega$ be a bounded Lipschitz domain such that the Poisson-Neumann problem \eqref{Neumann_problem} has $W^{2,p}$-regularity. There is a constant $C > 0$ such that
    \begin{equation} \vspace{-1pt}  
      \norma{\nabla z}{W^{1,p}(\Omega)} \leq C \norma{\Delta z}{L^p(\Omega)}, \ \forall z \in W^{2,p}(\Omega) \mbox{ such that } \partial_\eta z \Big |_{\Gamma} = 0.
      \label{obsv_norma_de_nabla_v}
    \vspace{-1pt} \end{equation} 
  \end{lemma}
  \begin{proof}[\bf Proof]
    Suppose that the result is false, that is, for each $n \in \mathbb{N}$ there is $z_n \in W^{2,p}(\Omega)$ with $\partial_\eta z_n \Big |_{\Gamma} = 0$ such that
    \begin{equation} \vspace{-1pt}  \label{negacao_resultado}
      \norma{\nabla z_n}{W^{1,p}(\Omega)} > n \norma{\Delta z_n}{L^p(\Omega)}.
    \vspace{-1pt} \end{equation} 
    Without loss of generality, we can take $z_n$ such that
    \begin{equation} \vspace{-1pt}  \label{sequencia_z_n}
      z^{\ast} = \int_{\Omega}{z_n \ dx} = 0 \mbox{ and } \norma{\nabla z_n}{W^{1,p}(\Omega)} = 1
    \vspace{-1pt} \end{equation} 
    Accounting for \eqref{negacao_resultado}, \eqref{sequencia_z_n} and Lemma \ref{lema_desigualdade_poincare_media_nula} we have $(z_n), (\nabla z_n)$ bounded in $W^{1,p}(\Omega)$ and
    \begin{equation} \vspace{-1pt}  \label{convergencia_Delta_z_n}
      \Delta z_n \longrightarrow \Delta z = 0 \mbox{ strongly in } L^p(\Omega).
    \vspace{-1pt} \end{equation} 
    Using the $W^{2,p}$-regularity of the Poisson-Neumann problem \eqref{Neumann_problem} we have
    \begin{equation*} \vspace{-1pt}
      \norma{z_n}{W^{2,p}(\Omega)} \leq C (\norma{\Delta z_n}{L^p(\Omega)} + \norma{z_n}{L^p(\Omega)})
    \vspace{-1pt} \end{equation*}  
    and thus $(z_n)$ is bounded in $W^{2,p}(\Omega)$. This allows us to conclude, using compactness results in Sobolev spaces, that there is $z \in W^{2,p}(\Omega)$ such that, up to a subsequence,
    \begin{equation} \vspace{-1pt}  \label{convergencia_z_n}
      z_n \longrightarrow z \mbox{ weakly in } W^{2,p}(\Omega) \mbox{ and strongly in } W^{1,p}(\Omega),
    \vspace{-1pt} \end{equation} 
    Using again the $W^{2,p}$-regularity of the Poisson-Neumann problem \eqref{Neumann_problem} we have
    \begin{equation*} \vspace{-1pt}
      \norma{z_i - z_j}{W^{2,p}(\Omega)} \leq C \norma{\Delta z_i - \Delta z_j}{L^p(\Omega)} + C \norma{z_i - z_j}{L^p(\Omega)}, \ \forall i,j \in \mathbb{N},
    \vspace{-1pt} \end{equation*}  
    and accounting for \eqref{convergencia_z_n} and \eqref{convergencia_Delta_z_n} we conclude that
    \begin{equation} \vspace{-1pt}  \label{convergencia_z_n_forte}
      z_n \longrightarrow z \mbox{ strongly in } W^{2,p}(\Omega).
    \vspace{-1pt} \end{equation} 
    Now, considering the properties of each $z_n$ and the convergences \eqref{convergencia_Delta_z_n} and \eqref{convergencia_z_n_forte} we have
    \begin{equation*} \vspace{-1pt}
      \Delta z = 0, \mbox{ with } \partial_\eta z \Big |_{\Gamma} = 0, \ z^{\ast} = 0 \mbox{ and } \norma{\nabla z}{W^{1,p}(\Omega)} = 1.
    \vspace{-1pt} \end{equation*}  
    But this is not possible because if $z$ satisfies $\Delta z = 0$, $\partial_\eta z \Big |_{\Gamma} = 0$ and $z^{\ast} = 0$, then we have $z \equiv 0$ and hence $\norma{\nabla z}{W^{1,p}(\Omega)} = 0$. Therefore we must have \eqref{obsv_norma_de_nabla_v}.
  \end{proof}
  
  \begin{hypothesis}
    For each $z \in H^2(\Omega)$ such that $\partial_{\eta} z \Big |_{\Gamma} = 0$, there is a sequence $\{ \rho_n \} \subset C^2(\overline{\Omega})$ such that $\partial_{\eta} \rho_n \Big |_{\Gamma} = 0$ and $\rho_n \to z$ in $H^2(\Omega)$.
    \label{hypothesis_density}
  \end{hypothesis}

\

  \begin{lemma}[\bf \cite{ViannaGuillen2023uniform}]
    Suppose that the Poisson-Neumann problem \eqref{Neumann_problem} has the $H^2$-regularity and assume that Hypothesis \ref{hypothesis_density} holds. Let $z \in H^2(\Omega)$ be such that $z \geq \alpha^2$, for some $\alpha^2 > 0$. Then there exist positive constants $C_1, C_2 > 0$, independent of $\alpha^2$, such that,
    \begin{align*}
      \int_{\Omega}{\norm{\Delta z}{}^2 \ dx} + \int_{\Omega}{\frac{\norm{\nabla z}{}^2}{z} \Delta z \ dx} & \geq C_1 \Big ( \int_{\Omega}{\norm{D^2 z}{}^2 \ dx} + \int_{\Omega}{\frac{\norm{\nabla z}{}^4}{z^2} \ dx} \Big ) - C_2 \int_{\Omega}{\norm{\nabla z}{}^2 \ dx}.
    \end{align*}
    \label{lema_termo_fonte_final}
  \end{lemma}  

  \begin{remark} \label{remark_hypotheses_domain}
    The hypotheses of Lemma \ref{lema_termo_fonte_final} related to the $H^2$-regularity of the Poisson-Neumann problem \eqref{Neumann_problem} and Hypothesis \ref{hypothesis_density} can be understood as hypotheses on the domain $\Omega$. Indeed, see \cite{Grisvard} and \cite[Appendix A.1]{ViannaGuillen2023uniform} for assumptions on $\Omega$ that imply $H^2$-regularity of the Poisson-Neumann problem \eqref{Neumann_problem} and the validity Hypothesis \ref{hypothesis_density}, respectively. Since the aforementioned hypotheses of Lemma \ref{lema_termo_fonte_final} are the most restrictive ones in this sense, from now on we assume that $\Omega$ is a bounded domain of $\mathbb{R}^3$ such that:
    \begin{enumerate}
      \item the Poisson-Neumann problem \eqref{Neumann_problem} has $H^2$-regularity and
      \item Hypothesis \ref{hypothesis_density} is valid.
    \end{enumerate}
  \end{remark}

\
  
  \begin{lemma}[\bf \cite{eyre1998unconditionally}]
    Let $z^n, z^{n-1} \in L^\infty(\Omega)$ and let $f: \mathbb{R} \rightarrow \mathbb{R}$ be a $C^2$ function. Then
    \begin{equation*} \vspace{-1pt}
      \int_{\Omega} \delta_t z^n \, f'(z^n) = \delta_t \int_{\Omega}{f(z^n) \ dx} + \frac1{2k} \int_{\Omega} f''(c^n(x))(z^n(x) - z^{n-1}(x))^2 \ dx,
    \vspace{-1pt} \end{equation*}  
    where $c^n(x)$ is an intermediate point between $z^n(x)$ and $z^{n-1}(x)$. In particular, if $f$ is convex then we have
    \begin{equation*} \vspace{-1pt}
     \int_{\Omega} \delta_t z^n\, f'(z^n)\, dx \geq \delta_t \int_{\Omega}{f(z^n) \, dx}.
    \vspace{-1pt} \end{equation*}  
    \label{lema_delta_t}
  \end{lemma}
  \begin{lemma}[\bf \cite{ViannaGuillen2023uniform}]
    Let $w_1$ and $w_2$ be non-negative real numbers. For each $s \geq 1$ we have
    \begin{equation*} \vspace{-1pt}
      \norm{w_2^s - w_1^s}{} \leq s \norm{w_2 + w_1}{}^{s-1} \norm{w_2 - w_1}{}.
    \vspace{-1pt} \end{equation*}  
    \label{lema_T^m_elevado_a_s}
  \end{lemma}
  Using Lemma \ref{lema_T^m_elevado_a_s}, we can prove the following.
  \begin{lemma}
     Let $p \in (1, \infty)$ and let $\{ w_n \}$ be a sequence of non-negative functions  such that $w_n \to w$ in $L^p(0,T;L^p(\Omega))$ as $n \to \infty$. Then, for every $r \in (1,p)$, $w_n^r \to w^r$ in $L^{p/r}(0,T;L^{p/r}(\Omega))$ as $n \to \infty$.
    \label{lema_convergencia_w_elevado_a_s}
  \end{lemma}
  
  The following result (Corollary $4$ of \cite{Simon1986compact}) establishes a criterion of compactness in Bochner spaces.
  
  \begin{lemma}[\bf Compactness in Bochner spaces]
    Let $X,B$ and $Y$ be Banach spaces, let
    \begin{equation*} \vspace{-1pt}
      F \subset \Big \{ f \in L^1(0,T;Y) \ \Big | \ \partial_t f \in L^1(0,T;Y) \Big \} \qquad \mbox{ and } \qquad \partial F/ \partial t = \Big \{ \partial_t f, \ \forall f \in F \Big \}.
    \vspace{-1pt} \end{equation*}  
    Suppose that $X \subset B \subset Y$, with compact embedding $X \subset B$ and continuous embedding $B \subset Y$. Let the set $F$ be bounded in $L^p(0,T;B) \cap L^1(0,T;X)$, for $1 < p \leq \infty$, and $\partial F/ \partial t$ be bounded in $L^1(0,T;Y)$. Then $F$ is relatively compact in $L^q(0,T;B)$, for $1 \leq q < p$.
    \label{lema_Simon}
  \end{lemma}



\section{Study of the \texorpdfstring{$\boldsymbol{(u,z)}$}{{\bf (u,z)}}-scheme \texorpdfstring{\eqref{NLTD}}{{\bf (u,z)}}}  \label{Sec:u-z}
  \subsection{Existence of solution of \texorpdfstring{\eqref{NLTD}}{{\bf (u,z)-scheme}}} 
  We remind that the solutions of \eqref{NLTD} are denoted by $(u^n,z^n)$ and recall Remark \ref{remark_dependencia_esquema_em_m}.
  
  \begin{theorem}{\bf (Existence of solution of \eqref{NLTD})}
    Suppose $(u^{n-1},z^{n-1}) \in L^2(\Omega) \times L^{\infty}(\Omega)$ with $u^{n-1}(x) \geq 0$ and $z^{n-1}(x) \geq \alpha$ $a.e.$ $x \in \Omega$. Then there is a solution $(u^n,z^n)$ of \eqref{NLTD} which satisfies $u^n(x) \geq 0$ and $\norma{z^{n-1}}{L^{\infty}(\Omega)} \geq z^n(x) \geq \alpha$ $a.e.$ $x \in \Omega$.
    \label{teo_existencia_NLTD}
  \end{theorem}
  \begin{proof}[\bf Proof]
    In order to avoid divisions by zero in some terms of \eqref{NLTD} and obtain $u^n(x) \geq 0$ $a.e.$ $x \in \Omega$, we define the lower truncation for $z$
    \begin{equation*} \vspace{-1pt}
      T_{\alpha}(z) = \left \{
      \begin{array}{rl}
        \alpha, & \mbox{if } z \leq \alpha, \\
        z, & \mbox{if } z \geq \alpha,
      \end{array}
      \right.
    \vspace{-1pt} \end{equation*}  
    and the lower-upper truncation for $u$
    \begin{equation*} \vspace{-1pt}
      T_0^m(u) = \left \{
      \begin{array}{rl}
        0, & \mbox{if } u \leq 0, \\
        u, & \mbox{if } u \in [0,m], \\
        m, & \mbox{if } u \geq m,
      \end{array}
      \right.
    \vspace{-1pt} \end{equation*}  
    Then, we consider the auxiliary problem
    \begin{equation} \vspace{-1pt}  
    \left\{
      \begin{array}{l}
        \delta_t u^n - \Delta u^n  =  \nabla \cdot \Big (  T_0^m(u^n) \nabla (z^n)^2 \Big ), \\
        \delta_t z^n - \dfrac{\norm{\nabla z^n}{}^2}{T_{\alpha}(z^n)} - \Delta z^n  = - \dfrac{1}{2} T_0^m(u^n)^s \left(z^n - \dfrac{\alpha^2}{T_{\alpha}(z^n)} \right),
      \end{array}
      \right.
      \label{NLTD_aux}
    \vspace{-1pt} \end{equation} 
    with the same boundary and initial conditions of \eqref{NLTD}.
    
    We prove the existence of a solution $(u^n,z^n)$ to \eqref{NLTD_aux} via Leray-Schauder fixed point theorem \cite{GilbargTrudinger}. Along this proof, we also have that any solution $(u^n,z^n)$ of \eqref{NLTD_aux} satisfies $u^n(x) \geq 0$ and $\norma{z^{n-1}}{L^{\infty}(\Omega)} \geq z^n(x) \geq \alpha$ $a.e.$ $x \in \Omega$, which implies that $T_{\alpha}(z^n) = z^n$, $T_0^m(u^n) = T^m(u^n)$ and therefore we conclude that $(u^n,z^n)$ is also a solution of \eqref{NLTD}. Now we proceed with the proof of existence for \eqref{NLTD_aux} which is divided in three steps.
    
  \
    
    \noindent {\bf Step 1 (Definition of the compact mapping $\boldsymbol{S}$):} For all $(\overline{u},\overline{z}) \in  W^{1,4}(\Omega) ^2$, we define $(u,z) = S(\overline{u},\overline{z}) \in H^2(\Omega) ^2$ as the solution of
    \begin{equation} \vspace{-1pt}  \label{NLTD_aux_Leray_Schauder-u}
        \dfrac{u}{k} - \Delta u = 2  \nabla \cdot \Big (  T_0^m(\overline{u})\, \overline{z}\, \nabla z \Big ) + \dfrac{u^{n-1}}{k}, 
        \vspace{-1pt} \end{equation} 
         \begin{equation} \vspace{-1pt}  \label{NLTD_aux_Leray_Schauder-z}
        \dfrac{z}{k} - \Delta z + \dfrac{1}{2} T_0^m(\overline{u})^s \left(z - \dfrac{\alpha^2}{T_{\alpha}(\overline{z})} \right) = \dfrac{\norm{\nabla \overline{z}}{}^2}{T_{\alpha}(\overline{z})} + \dfrac{z^{n-1}}{k}.
    \vspace{-1pt} \end{equation} 
    We can use standard results on linear elliptic problems to conclude that $(u,z) = S(\overline{u},\overline{z})$ is well defined. In fact, given $(\overline{u},\overline{z})$, we begin by solving the $z$-equation \eqref{NLTD_aux_Leray_Schauder-z}. Since $0 < \frac{1}{k} + \frac{1}{2} T_0^m(\overline{u})^s \leq \frac{1}{k} + \frac{m^s}{2}$, we first prove the existence of a weak solution $z \in H^1(\Omega)$ by means of the Lax-Milgram Theorem and then we use the $H^2$-regularity of the Poisson-Neumann problem \eqref{Neumann_problem} to prove that $z \in H^2(\Omega)$. Once the existence of $z \in H^2(\Omega)$ is proved, we have $\nabla \cdot \Big (  T_0^m(\overline{u})\, \overline{z} \,\nabla z \Big ) \in L^2(\Omega)$ and therefore we are able to solve the $u$-equation \eqref{NLTD_aux_Leray_Schauder-u}. Using again the $H^2$-regularity of the Poisson-Neumann problem \eqref{Neumann_problem}, we obtain $u \in H^2(\Omega)$. Hence, $S(\overline{u},\overline{z}) \in H^2(\Omega)^2$ and therefore $S$ is a compact mapping defined in $W^{1,4}(\Omega) ^2$.
    
  \
    
    \noindent {\bf Step 2 (Pointwise bounds for $\boldsymbol{u}$ and $\boldsymbol{z/\lambda}$ for any $\boldsymbol{(u,z) = \lambda S(u,z)}$):}
    Let $\lambda \in [0,1]$. We will study the pairs $(u,z)$ such that $(u,z) = \lambda S(u,z)$. If we consider $\lambda = 0$ then $S(u,z) \in H^2(\Omega)^2$ is well defined and $(u,z) = (0,0)$. Once the case where $\lambda = 0$ is treated, we consider $\lambda \in (0,1]$, therefore we can write $S(u,z) = (1/\lambda)^*(u,z)$ and we have $(u,z)$ satisfying
    \begin{equation} \vspace{-1pt}  
    \left\{
      \begin{array}{l}
        \dfrac{u}{k} - \Delta u =  \nabla \cdot \Big ( T_0^m(u) \nabla (z)^2 \Big ) + \lambda \dfrac{u^{n-1}}{k}, \\[12pt]
        \dfrac{1}{k} \dfrac{z}{\lambda} - \dfrac{\norm{\nabla z}{}^2}{T_{\alpha}(z)} - \Delta \dfrac{z}{\lambda}  = - \dfrac{1}{2} T_0^m(u)^s \left(\dfrac{z}{\lambda} - \dfrac{\alpha^2}{T_{\alpha}(z)} \right) + \dfrac{1}{k} z^{n-1}.
      \end{array}
      \right.
      \label{NLTD_aux_Leray_Schauder_lambda}
    \vspace{-1pt} \end{equation} 
    If we test the $u$-equation of \eqref{NLTD_aux_Leray_Schauder_lambda} by the negative part of $u$ defined, as $u_-(x) = \min{ \{ 0, u(x) \} }$, we conclude that $u \geq 0$. Now let $c = \norma{z^{n-1}}{L^{\infty}(\Omega)}$. We rewrite the $z$-equation of \eqref{NLTD_aux_Leray_Schauder_lambda} as
    \begin{equation} \vspace{-1pt}  
      \begin{array}{rl}
        & \dfrac{1}{k} (\dfrac{z}{\lambda} - c) - \lambda^2 \dfrac{\norm{\nabla (\dfrac{z}{\lambda} - c)}{}^2}{T_{\alpha}(z)} - \Delta (\dfrac{z}{\lambda} - c) = - \dfrac{T_0^m(u)^s}{2} (\dfrac{z}{\lambda} - c) + \dfrac{T_0^m(u)^s}{2} (\dfrac{\alpha^2}{T_{\alpha}(z)} - c) + \dfrac{(z^{n-1} - c)}{k}.
      \end{array}
      \label{NLTD_aux_Leray_Schauder_lambda_eq_z_menos_c}
    \vspace{-1pt} \end{equation} 
    Continuing, we test \eqref{NLTD_aux_Leray_Schauder_lambda_eq_z_menos_c} by the positive part of $(\dfrac{z}{\lambda} - c)$, defined as $(\dfrac{z}{\lambda} - c)_+(x) = \max{\{ 0, (\dfrac{z}{\lambda} - c)(x) \}} \geq 0$. Since $c = \norma{z^{n-1}}{L^{\infty}(\Omega)}$, we have $z^{n-1} - c \leq 0$. Moreover, note that $(\dfrac{z}{\lambda} - c)_+ \geq 0$ and $(\dfrac{z}{\lambda} - c)_+ \neq 0 \Longleftrightarrow z > \lambda c$. Then, an analysis taking account of the possible cases $\alpha < \lambda c$ and $\alpha \geq \lambda c$ leads us to
    \begin{equation*} \vspace{-1pt}
      \frac{(\dfrac{z}{\lambda} - c)_+}{T_{\alpha}(z)} \leq \dfrac{1}{\lambda}.
    \vspace{-1pt} \end{equation*}  
    Hence, reminding that $\lambda \in (0,1]$ and $T_{\alpha}(z) \geq \alpha$, if we test \eqref{NLTD_aux_Leray_Schauder_lambda_eq_z_menos_c} by $(z/\lambda - c)_+$ we obtain
    \begin{align*}
      \frac{1}{k} \norma{(\dfrac{z}{\lambda} - c)_+}{L^2(\Omega)}^2 + \norma{\nabla (\dfrac{z}{\lambda} - c)_+}{L^2(\Omega)}^2 \leq \lambda \norma{\nabla (\dfrac{z}{\lambda} - c)_+}{L^2(\Omega)}^2 + \frac{1}{2} T_0^m(u)^s (\alpha - c) (\dfrac{z}{\lambda} - c)_+.
    \end{align*}
    By hypothesis we have $z^{n-1} \geq \alpha$, which implies in particular $\alpha \leq c$ and then that $\dfrac{z}{\lambda} \leq \norma{z^{n-1}}{L^{\infty}(\Omega)}$.
    
    Next we prove an inferior bound for $z/\lambda$. Considering the definition of $T_{\alpha}$ and $\dfrac{\norm{\nabla z}{}^2}{T_{\alpha}(z)} + \dfrac{1}{k} z^{n-1} \geq \dfrac{\alpha}{k}$, which comes from the hypotheses of the theorem, we can use the $z$-equation of \eqref{NLTD_aux_Leray_Schauder_lambda} to write
    \begin{equation} \vspace{-1pt}  
      \dfrac{1}{k} (\dfrac{z}{\lambda} - \alpha) - \Delta (\dfrac{z}{\lambda} - \alpha) + \dfrac{1}{2} T_0^m(u^n)^s (\dfrac{z}{\lambda} - \alpha) \geq \dfrac{1}{2} T_0^m(u)^s (\dfrac{\alpha^2}{T_{\alpha}(z)} - \alpha).
      \label{NLTD_aux_z_ineq}
    \vspace{-1pt} \end{equation} 
    Now we test \eqref{NLTD_aux_z_ineq} by the negative part of $(z - \lambda \alpha)$, $(z - \lambda \alpha)_-$. Note that $(z - \lambda \alpha)_- \leq 0$ and $(\dfrac{z}{\lambda} - \alpha)_- \neq 0 \Longrightarrow z < \lambda \alpha \leq \alpha$. Therefore, testing \eqref{NLTD_aux_z_ineq} by $(\dfrac{z}{\lambda} - \alpha)_-$ and reminding that $\lambda \in (0,1]$, we obtain
    \begin{align*}
      & \frac{1}{k} \norma{(\dfrac{z}{\lambda} - \alpha)_-}{L^2(\Omega)}^2 + \norma{\nabla (\dfrac{z}{\lambda} - \alpha)_-}{L^2(\Omega)}^2 + \frac{1}{2} \norma{T_0^m(u)^{s/2} (\dfrac{z}{\lambda} - \alpha)_-}{L^2(\Omega)}^2 \\
      & \leq \frac{1}{2} \int_{\Omega} T_0^m(u)^s (\dfrac{\alpha^2}{T_{\alpha}(z)} - \alpha) (\dfrac{z}{\lambda} - \alpha)_- \ dx = \frac{1}{2} \int_{\Omega} T_0^m(u)^s (\dfrac{\alpha^2}{\alpha} - \alpha) (\dfrac{z}{\lambda} - \alpha)_- \ dx
    \end{align*}
    Thus we conclude that $(\dfrac{z}{\lambda} - \alpha)_- = 0$, that is, $\dfrac{z}{\lambda} \geq \alpha > 0$.

\

    \noindent {\bf Step 3 ($\boldsymbol{\lambda}$-independent bounds for any $\boldsymbol{(u,z) = \lambda S(u,z)}$):} As we mentioned before, we consider $\lambda \in (0,1]$ because if $\lambda = 0$ then we have $(u,z) = (0,0)$. Because of the upper bound for $z$ that we proved in the anterior step we have
    \begin{equation} \vspace{-1pt}  
      \frac{z}{\lambda} \mbox{ is bounded in } L^{\infty}(\Omega),
      \label{limitacao_z_L_infty_independente_lambda}
    \vspace{-1pt} \end{equation} 
    independently of $\lambda$. Then we can integrate the $z$-equation of \eqref{NLTD_aux_Leray_Schauder_lambda} and conclude that
    \begin{equation} \vspace{-1pt}  
      \norma{\nabla z}{L^2(\Omega)}^2 \leq C(k,\alpha^2,m,\norma{z^{n-1}}{L^{\infty}(\Omega)}).
      \label{estimativa_nabla_z_dependente_de_lambda}
    \vspace{-1pt} \end{equation} 
    Since $z > 0$, we can multiply the $z$-equation of \eqref{NLTD_aux_Leray_Schauder_lambda} by $\lambda\, T_{\alpha}(z)/z$. This gives us
    \begin{equation} \vspace{-1pt}  
      \dfrac{T_{\alpha}(z)}{k} - \lambda \dfrac{\norm{\nabla z}{}^2}{z} - \dfrac{T_{\alpha}(z)}{z} \Delta z = - \dfrac{1}{2} T_0^m(u)^s (T_{\alpha}(z) - \dfrac{\alpha^2}{z/\lambda}) + \dfrac{T_{\alpha}(z)}{z/\lambda} \dfrac{z^{n-1}}{k}.
      \label{NLTD_aux_Leray_Schauder_eq_z}
    \vspace{-1pt} \end{equation} 
    Now we test \eqref{NLTD_aux_Leray_Schauder_eq_z} by $- \Delta z$ and, using that $T_{\alpha}(z)/z \geq 1$, $z/\lambda \geq \alpha$ and \eqref{limitacao_z_L_infty_independente_lambda}, we obtain
    \begin{align*}
      & \frac{1}{k} \int_{\Omega} (T_{\alpha})'(z) \norm{\nabla z}{}^2 \ dx + (1 - \lambda) \int_{\Omega} \norm{\Delta z}{}^2 \ dx + \lambda \int_{\Omega} \norm{\Delta z}{}^2 \ dx + \lambda \int_{\Omega} \frac{\norm{\nabla z}{}^2}{z} \Delta z \ dx \\
      & \leq \frac{1}{2} \int_{\Omega} \norm{T_0^m(u)}{}^s \norm{T_{\alpha}(z) + \frac{\alpha^2}{\alpha}}{} \norm{\Delta z}{} \ dx + \int_{\Omega} \norm{\frac{T_{\alpha}(z)}{\alpha \ k} z^{n-1}}{} \norm{\Delta z}{} \ dx \\
      & \leq \dfrac{1}{\delta} C(m,\alpha^2,\norma{z^{n-1}}{L^{\infty}(\Omega)}) + \delta \norma{\Delta z}{L^2(\Omega)}^2.
    \end{align*}
    Then, applying Lemma \ref{lema_termo_fonte_final} and using \eqref{estimativa_nabla_z_dependente_de_lambda} and the $H^2$-regularity we obtain
    \begin{align*}
      & \frac{1}{k} \int_{\Omega} (T_{\alpha})'(z) \norm{\nabla z}{}^2 \ dx + (1 - \lambda) \int_{\Omega} \norm{\Delta z}{}^2 \ dx + C_1 \lambda ( \int_{\Omega} \norm{\Delta z}{}^2 \ dx + \int_{\Omega} \frac{\norm{\nabla z}{}^4}{z^2} \ dx ) \\
      & \leq C_2 \lambda \norma{\nabla z}{L^2(\Omega)}^2 + \dfrac{1}{\delta} C(m,\alpha^2,\norma{z^{n-1}}{L^{\infty(\Omega)}}) + \delta \norma{\Delta z}{L^2(\Omega)}^2 \\
      & \leq \dfrac{1}{\delta} C(k,m,\alpha^2,\norma{z^{n-1}}{L^{\infty(\Omega)}}) + \delta \norma{\Delta z}{L^2(\Omega)}^2.
    \end{align*}
    Taking into account that $\D{\int_{\Omega}} (T_{\alpha})'(z) \norm{\nabla z}{}^2 \ dx \geq 0$ and denoting $C = \min{\{ 1, C_1 \} }$ we have $1 + (C_1 - 1)\lambda \geq C$ for all $\lambda \in [0,1]$ and
    \begin{align*}
      & C \int_{\Omega} \norm{\Delta z}{}^2 \ dx + C \lambda \int_{\Omega} \frac{\norm{\nabla z}{}^4}{z^2} \ dx \leq \dfrac{1}{\delta} C(m,\alpha^2,\norma{z^{n-1}}{L^{\infty(\Omega)}}) + \delta \norma{\Delta z}{L^2(\Omega)}^2
    \end{align*}
    and therefore we can choose $\delta > 0$ small enough such that
    \begin{equation} \vspace{-1pt}  
      \Delta z \mbox{ is bounded in } L^2(\Omega).
      \label{limitacao_Delta_z_independente_lambda}
    \vspace{-1pt} \end{equation} 
    Because of the homogeneous Neumann boundary conditions, the $H^2$-regularity and the $\lambda$-uniform bounds \eqref{limitacao_z_L_infty_independente_lambda} and \eqref{limitacao_Delta_z_independente_lambda}, we conclude that
    \begin{equation} \vspace{-1pt}  
      z \mbox{ is bounded in } H^2(\Omega).
      \label{limitacao_z_H_2}
    \vspace{-1pt} \end{equation} 
    By Sobolev inequality in $3D$ domains, the latter implies that
    \begin{equation} \vspace{-1pt}  
      z \mbox{ is bounded in } W^{1,4}(\Omega).
      \label{limitacao_z_W_1_4}
    \vspace{-1pt} \end{equation} 
    Considering the $\lambda$-independent bounds $0 \leq z \leq \norma{z^{n-1}}{L^{\infty}(\Omega)}$, $T_0^m(u) \leq m$ and \eqref{limitacao_z_H_2}, we can test the $u$-equation \eqref{NLTD_aux_Leray_Schauder-u} by $u$ and prove that $u \mbox{ is bounded in } H^1(\Omega).$ Using again the $\lambda$-independent bounds $0 \leq z \leq \norma{z^{n-1}}{L^{\infty}(\Omega)}$, $T_0^m(u) \leq m$, \eqref{limitacao_z_H_2} and \eqref{limitacao_z_W_1_4}, we have that the chemotaxis term $2  \nabla \cdot \Big (  T_0^m(\overline{u}) \overline{z} \nabla z \Big )$ is bounded in $L^2(\Omega)$. Then we can test the $u$-equation \eqref{NLTD_aux_Leray_Schauder-u} by $- \Delta u$, obtaining $u$ is bounded in $H^2(\Omega)$, which implies, in particular, that
    \begin{equation} \vspace{-1pt}  
      u \mbox{ is bounded in } W^{1,4}(\Omega).
      \label{limitacao_u_W_1_4}
    \vspace{-1pt} \end{equation} 
    With \eqref{limitacao_z_W_1_4} and \eqref{limitacao_u_W_1_4} we finally conclude, using the Leray-Schauder fixed point theorem \cite{GilbargTrudinger} that the auxiliary problem \eqref{NLTD_aux} has a solution $(u,z)$. Because of the properties showed along the steps of the proof we also conclude that $(u,z) \in H^2(\Omega)^2$, $u(x) \geq 0$ and $\norma{z^{n-1}}{L^{\infty}(\Omega)} \geq z(x) \geq \alpha$ $a.e.$ $x \in \Omega$. Therefore we have $T_{\alpha}(z) = z$ and $T_0^m(u) = T^m(u)$, which implies that the solution $(u,z)$ of \eqref{NLTD_aux} is a solution of \eqref{NLTD}, finishing the proof of existence of solution.
  \end{proof}
  

\subsection{First uniform in time estimates}
  
  The following direct estimates and the energy inequalities obtained in this subsection are valid for any solution $(u^n,z^n)$ of \eqref{NLTD} given by Theorem \ref{teo_existencia_NLTD}.
  
  \begin{lemma}{\bf ( $\boldsymbol{(m,k,n)}$-uniform estimates)}
    Let $(u^n,z^n)$ be a solution of \eqref{NLTD}. Then we have
    \begin{enumerate}
      \item $\D{\int_{\Omega}} u^n \ dx = \D{\int_{\Omega}} u^0 \ dx$, for all $n \in \mathbb{N}$;
      \item $\norma{z^n}{L^2(\Omega)}^2 + \D{\sum_{j=1}^n} \norma{z^j - z^{j-1}}{L^2(\Omega)}^2 \leq \norma{z^0}{L^2(\Omega)}^2$, for all $n \in \mathbb{N}$;
      \item $k \D{\sum_{j=1}^n}{\norma{\nabla z^j}{L^2(\Omega)}^2} \leq \frac{1}{4 \alpha^2} \norma{v^0 + \alpha^2}{L^2(\Omega)}^2$, for all $n \in \mathbb{N}$.
    \end{enumerate}
    \label{lema_limitacao_grad_z_uniforme_wrt_k_m_T}
  \end{lemma}
  \begin{proof}[\bf Proof]
    The proof of $1$ is achieved by integrating the $u^n$-equation of \eqref{NLTD}.
    
    For the items $2$ and $3$ we take the product of the $z^n$-equation of \eqref{NLTD} by $z^n$. We obtain
    \begin{equation} \vspace{-1pt}  
      \delta_t (z^n)^2 + \frac{1}{k} (z^n - z^{n-1})^2 - \Delta (z^n)^2 + T^m(u^n)^s ((z^n)^2 - \alpha^2) = 0.
      \label{eq_lema_limitacao_grad_z}
    \vspace{-1pt} \end{equation} 
    Since $((z^n)^2 - \alpha^2) \geq 0$ and $\D{\int_{\Omega}} \Delta (z^n)^2 \ dx = 0$, by integrating \eqref{eq_lema_limitacao_grad_z} we prove item $2$. On the other hand, testing \eqref{eq_lema_limitacao_grad_z} by $k (z^n)^2$ leads us to
    \begin{equation*} \vspace{-1pt}
      k \delta_t \norma{(z^n)^2}{L^2(\Omega)}^2 + k \norma{\nabla (z^n)^2}{L^2(\Omega)}^2 \leq 0.
    \vspace{-1pt} \end{equation*}  
    Then, summing up from $j = 1$ to $n$ gives us
    \begin{equation*} \vspace{-1pt}
      k \sum_{j=1}^n{\norma{\nabla (z^j)^2}{L^2(\Omega)}^2} \leq \norma{(z^0)^2}{L^2(\Omega)}^2 = \norma{v^0 + \alpha^2}{L^2(\Omega)}^2.
    \vspace{-1pt} \end{equation*}  
    Now using that $z^j \geq \alpha$ we have
    \begin{equation*} \vspace{-1pt}
      \int_{\Omega} \norm{\nabla z^j}{}^2 \ dx = \int_{\Omega} \frac{(z^j)^2}{(z^j)^2} \norm{\nabla z^j}{}^2 \ dx \leq \int_{\Omega} \frac{1}{4 \alpha^2} \norm{\nabla (z^j)^2}{}^2 \ dx,
    \vspace{-1pt} \end{equation*}  
    hence we obtain $3$.
  \end{proof}
    

  \subsection{Energy inequality}
    
    Now we turn to the energy inequalities giving for $s \in [1,2)$ in Lemma \ref{lemma_estimativa_u_v_m_1_s_intermediario} and $s \geq 2$ in Lemma \ref{lemma_estimativa_u_v_m_1_s_geq_2}. We will need the following lemma.
    
    \begin{lemma}
      Any solution $(u^n,z^n)$ of \eqref{NLTD}, satisfies the inequality
      \begin{equation} \vspace{-1pt}  
        \begin{array}{c}
          \dfrac{1}{2} \delta_t \norma{\nabla z^n}{L^2(\Omega)}^2 + \dfrac{1}{2 k} \norma{\nabla z^n - \nabla z^{n-1}}{L^2(\Omega)}^2 + C_1 \Big ( \D{\int_{\Omega}}{\norm{D^2 z^n}{}^2 \ dx} + \D{\int_{\Omega}}{\frac{\norm{\nabla z^n}{}^4}{(z^n)^2} \ dx} \Big ) \\ + \dfrac{1}{2} \D{\int_{\Omega}{T^m(u^n)^s \norm{\nabla z^n}{}^2 \ dx}} \leq \dfrac{s}{4} \D{\int_{\Omega}}{T^m(u^n)^{s-1} \nabla (z^n)^2 \cdot \nabla T^m(u^n) \ dx} \\
          + \dfrac{s \alpha}{2} \D{\int_{\Omega}}{T^m(u^n)^{s-1} \norm{\nabla z^n}{} \norm{\nabla T^m(u^n)}{} \ dx} + C_2 \D{\int_{\Omega}}{\norm{\nabla z^n}{}^2 \ dx}.
        \end{array}
        \label{eq_tratamento_equacao_z}
      \vspace{-1pt} \end{equation} 
      \label{lema_tratamento_equacao_z}
    \end{lemma}
    \begin{proof}[\bf Proof]
      We begin by testing the $z^n$-equation of \eqref{NLTD} by $- \Delta z^n$. This gives us
      \begin{align*}
        & \frac{1}{2} \delta_t \norma{\nabla z^n}{L^2(\Omega)}^2 + \dfrac{1}{2 k} \norma{\nabla z^n - \nabla z^{n-1}}{L^2(\Omega)}^2 + \norma{\Delta z^n}{L^2(\Omega)}^2 \\
        & \qquad + \int_{\Omega}{\frac{\norm{\nabla z^n}{}^2}{z^n} \Delta z^n \ dx} + \frac{1}{2} \int_{\Omega}{(1 + \frac{\alpha^2}{(z^n)^2}) T^m(u^n)^s \norm{\nabla z^n}{}^2 \ dx} \\
        & = \frac{s}{4} \int_{\Omega}{T^m(u^n)^{s-1} \nabla (z^n)^2 \cdot \nabla T^m(u^n) \ dx} + \dfrac{s}{2} \alpha^2 \int_{\Omega}{\frac{T^m(u^n)^{s-1}}{z^n} \nabla z^n \cdot \nabla T^m(u^n) \ dx}.
      \end{align*}
      Then, estimating the last term on the left hand side by bellow, the last term on the right hand side by above and applying Lemma \ref{lema_termo_fonte_final}, we arrive at the desired inequality.
    \end{proof}
    
    Next we will obtain a local energy inequality for $(u^n, z^n)$, first for $s \in [1,2)$ and then for $s \geq 2$. We consider the function $f_m$ defined by $f_m(r) = \D{\int_0^r}{f'_m(\theta) \ d\theta}$, where
      \begin{equation*} \vspace{-1pt}
        f_m'(r) = \left \{ 
        \begin{array}{rl}
          ln(T^m(r)), & \mbox{if } s = 1, \mbox{ for } r > 0, \\
          \dfrac{T^m(r)^{s-1}}{(s-1)}, & \mbox{if } s > 1. 
        \end{array}
        \right.
      \vspace{-1pt} \end{equation*}  
    
    \begin{lemma}[\bf Energy inequality for $\boldsymbol{s \in [1,2)}$]
      Any solution $(u^n, z^n)$ of the problem \eqref{NLTD} satisfies, for sufficiently small $\alpha^2 > 0$,
      \begin{equation} \vspace{-1pt}  
        \begin{array}{c}
          \delta_t \Big [ \dfrac{s}{4} \D{\int_{\Omega}}{f_m(u^n) \ dx} + \frac{1}{2} \norma{\nabla z^n}{L^2(\Omega)}^2 \Big ] + \dfrac{1}{2 k} \norma{\nabla z^n - \nabla z^{n-1}}{L^2(\Omega)}^2 \\
          + \dfrac{1}{4} \D{\int_{\Omega}}{T^m(u^n)^s \norm{\nabla z^n}{}^2 \ dx} + C_1 \Big ( \D{\int_{\Omega}}{\norm{D^2 z^n}{}^2 \ dx} + \D{\int_{\Omega}}{\frac{\norm{\nabla z^n}{}^4}{(z^n)^2} \ dx} \Big ) \leq C \norma{\nabla z^n}{L^2(\Omega)}^2.
        \end{array}
        \label{estimativa_u_v_m_1_s_intermediario}
      \vspace{-1pt} \end{equation} 
      \label{lemma_estimativa_u_v_m_1_s_intermediario}
    \end{lemma}
    \begin{proof}[\bf Proof]
      The proof follows the same ideas of the analogous result that was proved in \cite{ViannaGuillen2023uniform} for the truncated model. We are going to show the main steps of the proof, calling the attention to the differences that appear due to the fact that we are dealing with the time discrete scheme \eqref{NLTD}. In fact, we begin by considering the sequence $\{ 1/j \}_{j \in \mathbb{N}}$ and by testing the $u^n$-equation of \eqref{NLTD} by the function $f_{m,j}'(u^n) = f_m'(u^n + 1/j)$. A difference appears in the treatment of the term of the discrete time derivative. Given that $f_{m,j}$ is convex, we use Lemma \ref{lema_delta_t} to obtain
      \begin{align*}
        \delta_t \int_{\Omega}{f_{m,j}(u^n) \ dx} + \int_{\Omega}{\dfrac{(T^m)'(u^n)}{(T^m(u^n) + 1/j)^{2-s}} \norm{\nabla u^n}{}^2 \ dx} \leq \prodl{\dfrac{T^m(u^n)}{(T^m(u^n) + 1/j)^{2-s}} \nabla (z^n)^2}{\nabla T^m(u^n)}.
      \end{align*}
      Then we can follow the ideas in \cite{ViannaGuillen2023uniform} until the point that we reach the inequality
      \begin{equation} \vspace{-1pt}  \label{eq_estimativa_j}
        \begin{array}{l}
          \delta_t \Big [ \dfrac{s}{4} \D{\int_{\Omega}} f_{m,j}(u^n) \ dx + \dfrac{1}{2} \norma{\nabla z^n}{L^2(\Omega)}^2 \Big ] + \dfrac{1}{2 k} \norma{\nabla z^n - \nabla z^{n-1}}{L^2(\Omega)}^2 \\
          + C_1 \Big ( \D{\int_{\Omega}}{\norm{D^2 z^n}{}^2 \ dx} + \D{\int_{\Omega}}{\frac{\norm{\nabla z^n}{}^4}{(z^n)^2} \ dx} \Big ) + \dfrac{1}{2} \D{\int_{\Omega}}{\Big [ T^m(u^n)^s - \dfrac{1}{2} (T^m(u^n) + 1/j)^s \Big ] \norm{\nabla z^n}{}^2 \ dx} \\
          \leq C \norma{\nabla z^n}{L^2(\Omega)}^2 + \dfrac{s}{4} \D{\int_{\Omega}{\Big [ T^m(u^n)^{s-1} - (T^m(u^n) + 1/j)^{s-1} \Big ] \nabla T^m(u^n) \cdot \nabla (z^n)^2 \ dx}}.
        \end{array}
      \vspace{-1pt} \end{equation} 
      
      Finally we pass to the limit as $j \to \infty$ in \eqref{eq_estimativa_j}. The presence of the discrete time derivative $\delta_t$ instead of $\partial_t$ is what allows us to give a unified treatment for the case $s \in [1,2)$, differently from \cite{ViannaGuillen2023uniform}. We proceed with the passage to the limit term by term. We detail the passage to the limit in the term which involves the discrete time derivative,
      \begin{equation} \vspace{-1pt}  
        \delta_t \int_{\Omega}{f_{m,j}(u^n(x)) \ dx}.
        \label{termo_derivada}
      \vspace{-1pt} \end{equation} 
      We define the functions $g_{m,j}, g_m, G \in L^1(\Omega)$ by $g_{m,j}(x) = \delta_t f_{m,j}(u^n(x))$, $g_m(x) = \delta_t f_m(u^n(x))$ and $G(x) = \norm{g_{m,1}(x)}{}$. Then, for almost every $x \in \Omega$, $g_{m,j}(x) \to g_m(x)$ as $j \to \infty$ with $\norm{g_{m,j}(x)}{} \leq G(x)$ for all $j \in \mathbb{N}$. Therefore, using the Dominated Convergence Theorem, we conclude that
      \begin{align*}
         \lim_{j \to \infty}{\delta_t \int_{\Omega}{f_{m,j}(u^n(x)) \ dx}} & = \lim_{j \to \infty}{\int_{\Omega}{g_{m,j}(x) \ dx}} = \int_{\Omega}{g_m(x) \ dx} = \delta_t \int_{\Omega}{f_m(u^n(x)) \ dx}.
       \end{align*}
       
       For the other terms of \eqref{eq_estimativa_j}, one can again follow \cite{ViannaGuillen2023uniform}, take the limit as $j \to \infty$ and reach the desired result.
      \end{proof}
      
      \begin{lemma}[\bf Energy inequality for $\boldsymbol{s \geq 2}$] \label{lemma_estimativa_u_v_m_1_s_geq_2}
        The solution $(u^n, z^n)$ of the problem \eqref{NLTD} satisfies
        \begin{equation} \vspace{-1pt}  \label{estimativa_u_v_m_1_s_geq_2}
          \begin{array}{c}
            \delta_t \Big [ \dfrac{s}{4} \D{\int_{\Omega}} f_m(u^n) \ dx  + \frac{1}{2} \norma{\nabla z^n}{L^2(\Omega)}^2 \Big ] + \dfrac{1}{2 k} \norma{\nabla z^n - \nabla z^{n-1}}{L^2(\Omega)}^2 + \D{\int_{\Omega}}{\norm{\nabla [T^m(u^n)]^{s/2}}{}^2  dx} \\
            + \dfrac{1}{4} \D{\int_{\Omega}}{T^m(u^n)^s \norm{\nabla z^n}{}^2 \ dx} + C_1 \Big ( \D{\int_{\Omega}}{\norm{D^2 z^n}{}^2 \ dx} + \D{\int_{\Omega}}{\frac{\norm{\nabla z^n}{}^4}{(z^n)^2} \ dx} \Big ) \leq C \norma{\nabla z^n}{L^2(\Omega)}^2.
          \end{array}
        \vspace{-1pt} \end{equation} 
      \end{lemma}
      \begin{proof}[\bf Proof]
        The proof follows the same ideas of the analogous result that was proved in \cite{ViannaGuillen2023uniform} for the truncated model. Having in mind that we use Lemma \ref{lema_delta_t} to treat the term which involves the discrete time derivative $\delta_t$, we refer the reader to \cite{ViannaGuillen2023uniform} for the details of the proof.
      \end{proof}
    
    The energy inequalities \eqref{estimativa_u_v_m_1_s_intermediario} and \eqref{estimativa_u_v_m_1_s_geq_2} allow us to obtain $(m,k,n)$-independent estimates for the function $z^n$ in the next subsection.


  \section{Energy estimates and passage to the limit as \texorpdfstring{$\boldsymbol{(m,k) \to (\infty,0)}$}{(m,k) goes to (infinity,0)}}
    \label{subsec:passage_limit_m_to_infty}
    
    Now we use the global in time functions introduced in \eqref{funcao_u^kr_u^k}. With these functions we are able to rewrite \eqref{NLTD} as the following differential system, $a.e.$ in $(t,x) \in (0,\infty) \times \Omega$,
    \begin{equation} \vspace{-1pt}  
    \left\{
      \begin{array}{l}
        \partial_t \tilde{u}_m^k - \Delta u_m^k  =  \nabla \cdot \Big (  T^m(u_m^k) \nabla (z_m^k)^2 \Big ), \\
        \partial_t \tilde{z}_m^k - \dfrac{\norm{\nabla z_m^k}{}^2}{z_m^k} - \Delta z_m^k  = - \dfrac{1}{2} T^m(u_m^k)^s \left(z_m^k - \dfrac{\alpha^2}{z_m^k}\right).
      \end{array}
      \right.
      \label{NLTD_equiv}
    \vspace{-1pt} \end{equation}
    
    In this subsection we are going to prove $(m,k)$-independent estimates for the $(m,k)$-sequences $\{ u_m^k \}$, $\{ z_m^k \}$, $\{ \tilde{u}_m^k \}$ and $\{ \tilde{z}_m^k \}$, which will also be uniform in time.
    
    First, in Subsection \ref{subsec:estimativas_z_m}, we obtain estimates for $\nabla z_m^k$ from the energy inequalities \eqref{estimativa_u_v_m_1_s_intermediario} and \eqref{estimativa_u_v_m_1_s_geq_2}. Next we prove bounds for $\tilde{u}_m^k$ and $u_m^k$ and pass to the limit in \eqref{NLTD_equiv} as $(m,k) \to (\infty,0)$, considering the cases $s \in [1,2)$ and $s \geq 2$, separately.
    
    
    \subsection{Estimates for \texorpdfstring{$\boldsymbol{\nabla z_m^k}$}{grad(z\_m\^\{k,r\})}}
    \label{subsec:estimativas_z_m}
      
      Given $(u^n,z^n)$ a solution of \eqref{NLTD}, let us define the energy
      \begin{equation} \vspace{-1pt}  
        E_m^n = \frac{s}{4} \int_{\Omega}{f_m(u^n(t,x)) \ dx} + \frac{1}{2} \int_{\Omega}{\norm{\nabla z^n(t,x)}{}^2} \ dx.
        \label{energia_E_m}
      \vspace{-1pt} \end{equation} 
      
      We use the regularity of the initial data $u^0$ and $z^0$ in order to conclude that the initial energy $E_m^0$, is also bounded, independently of $(m,k)$. If we consider either \eqref{estimativa_u_v_m_1_s_intermediario} or \eqref{estimativa_u_v_m_1_s_geq_2}, multiply it by $k$ and sum from $j = 1$ to $n$ we obtain 
      \begin{equation} \vspace{-1pt}  
        \begin{array}{c}
          \dfrac{s}{4} \D{\int_{\Omega}} f_m(u^n) \ dx  + \frac{1}{2} \norma{\nabla z^n}{L^2(\Omega)}^2 + \dfrac{1}{2} \D{\sum_{j=1}^n}{\norma{\nabla z^j - \nabla z^{j-1}}{L^2(\Omega)}^2} \\
          \quad + k \D{\sum_{j=1}^n}{\D{\int_{\Omega}}{T^m(u^j)^s \norm{\nabla z^j}{}^2 \ dx}} + C_1 k \D{\sum_{j=1}^n}{\Big ( \D{\int_{\Omega}}{\norm{D^2 z^j}{}^2 \ dx} + \D{\int_{\Omega}}{\frac{\norm{\nabla z^j}{}^4}{(z^j)^2} \ dx} \Big )} \\
          \leq C k \D{\sum_{j=1}^n}{\norma{\nabla z^j}{L^2(\Omega)}^2} + \dfrac{s}{4} \D{\int_{\Omega}} f_m(u^0) \ dx  + \frac{1}{2} \norma{\nabla z^0}{L^2(\Omega)}^2,
        \end{array}
        \label{estimativa_nabla_z_m}
        \vspace{-1pt} \end{equation} 
      Thus, recalling that $z^j \geq \alpha$ in $\Omega$ and that Lemma \ref{lema_limitacao_grad_z_uniforme_wrt_k_m_T}.3 and \eqref{estimativa_nabla_z_m} are valid for any $n \in \mathbb{N}$, we conclude that
      \begin{equation} \vspace{-1pt}  \label{limitacao_aux_Dz_m_s_intermediario}
        \nabla z_m^k \mbox{ is bounded in } L^{\infty}(0,\infty;L^2(\Omega)) \cap L^4(0,\infty;L^4(\Omega)),
      \vspace{-1pt} \end{equation} 
      \begin{equation} \vspace{-1pt}  \label{limitacao_aux_Delta_z_m_s_intermediario}
        T^m(u_m^k)^{s/2} \nabla z_m^k \mbox{ and } \Delta z_m^k \mbox{ are bounded in } L^2(0,\infty;L^2(\Omega)),
      \vspace{-1pt} \end{equation} 
      \begin{equation} \vspace{-1pt}  \label{limitacao_diferenca_z_H_1}
        \D{\sum_{j=1}^{\infty}}{\norma{z^j - z^{j-1}}{H^1(\Omega)}^2} \leq C.
      \vspace{-1pt} \end{equation} 
      Using \eqref{limitacao_diferenca_z_H_1}, we prove the following.
      \begin{lemma}
        There is a positive constant $C$, independent of $m$ and $k$, such that
        \begin{equation} \vspace{-1pt}  \label{convergencia_diferenca_z_m_em_H_1}
          \norma{z_m^k - \tilde{z}_m^k}{L^2(0,\infty;H^1(\Omega))}^2 \leq C\, k.
        \vspace{-1pt} \end{equation} 
      \end{lemma}
      \begin{proof}[\bf Proof]
        From the definition of $z_m^k$ and $\tilde{z}_m^k$, we observe that $z_m^k(t) - \tilde{z}_m^k(t) = \dfrac{(t_n - t)}{k}(z^n - z^{n-1})$, for $t \in (t_{n-1},t_n)$. If $t \in (t_{j-1}, t_j)$ then $0 \leq t_j - t \leq k$ and hence
        \begin{align*}
          \norma{z_m^k - \tilde{z}_m^k}{L^2(0,\infty;H^1(\Omega))}^2 & = \sum_{j=1}^{\infty}{ \int_{t_{j-1}}^{t_j} \dfrac{(t_n - t)}{k} \norma{z^j - z^{j-1}}{H^1(\Omega)}^2} \ dt \leq k \sum_{j=1}^{\infty}{\norma{z^j - z^{j-1}}{H^1(\Omega)}^2}.
        \end{align*}
        Therefore, using \eqref{limitacao_diferenca_z_H_1} we conclude the proof.
      \end{proof}
      
      In particular, since $z_m^k(t,\cdot) \in H^2(\Omega)$ and $\dfrac{\partial}{\partial \eta} z^n(t,\cdot) |_{\Gamma} = 0$, it stems from \eqref{limitacao_aux_Delta_z_m_s_intermediario}, the $H^2$-regularity of the Poisson-Neumann problem \eqref{Neumann_problem} and \eqref{obsv_norma_de_nabla_v} that
      \begin{equation} \vspace{-1pt}  
        \nabla z_m^k \mbox{ is bounded in } L^2(0,\infty;H^1(\Omega)).
        \label{limitacao_aux_Dz_m_s_intermediario_2}
      \vspace{-1pt} \end{equation} 
      
      Using the results obtained until this point we analyze the existence of solutions of \eqref{problema_P}, first for $s \in [1,2)$ and then for $s \geq 2$.
      
      
      
      \subsection{Estimates for \texorpdfstring{$\boldsymbol{(u_m^k,z_m^k)}$}{(u\_m\^\{k,r\},z\_m\^\{k,r\})} and passage to the limit for \texorpdfstring{$\boldsymbol{s \in [1,2)}$}{s in [1,2)}}
      
      \label{subsec:estimativas_u_m_s<2}
      
      Let
      \begin{equation*} \vspace{-1pt}
        \qquad f'(r) =  \left \{
        \begin{array}{cc}
          ln(r) & \mbox{ if } s = 1, \ \forall r > 0, \\
          r^{s-1}/(s-1) & \mbox{ if } s \in (1,2),
        \end{array}
        \right .
      \vspace{-1pt} \end{equation*}  
      \begin{equation*} \vspace{-1pt}
        f(r) =  \int_0^r{f'(\theta) \ d\theta} = \left \{
        \begin{array}{cc}
          r ln(r) - r & \mbox{ if } s = 1, \ \forall r > 0, \\
          r^s/s(s-1) & \mbox{ if } s \in (1,2).
        \end{array}
        \right .
      \vspace{-1pt} \end{equation*}  
      Notice that $f''(r) = r^{s-2}, \ \forall r > 0$, in all cases.
      
      We test the $u^n$-equation of \eqref{NLTD} by $f'(u^n + 1)$. Using Lemma \ref{lema_delta_t} we obtain
      \begin{align*}
        & \delta_t \int_{\Omega} f(u^n + 1) \ dx + \frac{1}{2 k} \int_{\Omega} f''(c^n) (u^n - u^{n-1})^2 \ dx + \dfrac{4}{s^2} \int_{\Omega}{\norm{\nabla [u^n + 1]^{s/2}}{}^2 \ dx} + \\
        & = 2 \int_{\Omega}{T^m(u^n) (u^n + 1)^{s/2-1} z^n \nabla z^n  \cdot \nabla u^n \,  (u^n + 1)^{s/2-1} \ dx} \\
        & \leq \dfrac{4}{s} \int_{\Omega}{\dfrac{T^m(u^n)^{1-s/2}}{(u^n + 1)^{1-s/2}} T^m(u^n)^{s/2} z^n \nabla z^n \cdot \nabla [u^n + 1]^{s/2} \ dx} \\
        & \leq \dfrac{4}{s} \norma{z^0}{L^{\infty}(\Omega)} \Big ( \int_{\Omega}{T^m(u^n)^s \norm{\nabla z^n}{}^2 \ dx} \Big )^{1/2} \Big ( \int_{\Omega}{\norm{\nabla [u^n + 1]^{s/2}}{}^2 \ dx} \Big )^{1/2}
      \end{align*}
      and thus we have
      \begin{align*}
        & \delta_t \int_{\Omega}{f(u^n + 1) \ dx} + \frac{1}{2 k} \int_{\Omega} f''(c^n) (u^n - u^{n-1})^2 \ dx \\
        & + C_1 \int_{\Omega}{\norm{\nabla [u^n + 1]^{s/2}}{}^2 \ dx} \leq C_2 \int_{\Omega}{T^m(u^n)^s \norm{\nabla z^n}{}^2 \ dx}.
      \end{align*}
      
      Multiplying by $k$ and summing up from $0$ to $n$, for any $n \in \mathbb{N}$, we obtain
      \begin{align*}
        & \int_{\Omega}{f(u^n + 1) \ dx} + \frac{1}{2} \sum_{j=1}^n \int_{\Omega} f''(c^j) (u^j - u^{j-1})^2 \ dx + k \sum_{j=1}^n \dfrac{2}{s^2} \int_{\Omega}{\norm{\nabla [u^j + 1]^{s/2}}{}^2 \ dx} \\
        & \leq C k \sum_{j=1}^n \int_{\Omega}{T^m(u^j)^s \norm{\nabla z^j}{}^2 \ dx} + \int_{\Omega}{f(u^0 + 1) \ dx}.
      \end{align*}
      
      Then, because of \eqref{limitacao_aux_Delta_z_m_s_intermediario} and the definitions of $f$, $u_m^k$ and $\tilde{u}_m^k$ we conclude that
      \begin{equation} \vspace{-1pt}  
        \sum_{j=1}^{\infty} \int_{\Omega} f''(c^j) (u^j - u^{j-1})^2 \ dx \leq C,
        \label{limitacao_diferenca_u^n_u^n-1_s<2}
      \vspace{-1pt} \end{equation} 
      \begin{equation} \vspace{-1pt}  
          (u_m^k + 1)^{s/2} \mbox{ are bounded in } L^{\infty}(0,\infty;L^2(\Omega)),
          \label{limitacao_u_m_L^s}
      \vspace{-1pt} \end{equation} 
      in particular,
      \begin{equation} \vspace{-1pt}  
        u_m^k \mbox{ are bounded in } L^{\infty}(0,\infty;L^s(\Omega))
        \label{limitacao_u_m_L^s_2},
      \vspace{-1pt} \end{equation} 
      \begin{equation} \vspace{-1pt}  
        \nabla [u_m^k + 1]^{s/2} \mbox{ is bounded in } L^2(0,\infty;L^2(\Omega)).
        \label{limitacao_nabla_u_m_elevado_a_s_sobre_2}
      \vspace{-1pt} \end{equation} 
      
      Using these bounds we prove the following.
      \begin{lemma} \label{lema_f''}
        There is a positive constant $C$, independent of $m$ and $k$, such that
        \begin{equation} \vspace{-1pt}  
          \norma{u_m^k - \tilde{u}_m^k}{L^2(0,\infty;L^s(\Omega))}^2 \leq C\, k.
          \label{convergencia_diferenca_u^n_u^n-1_s<2}
        \vspace{-1pt} \end{equation} 
      \end{lemma}
      \begin{proof}[\bf Proof]
      
      \
      
        {\bf Step 1:} We remind that in \eqref{limitacao_diferenca_u^n_u^n-1_s<2} we have $f''(c^j) = (c^j)^{s-2}$, where $s \in [1,2)$ and, for each $j$ and for each $x$, $c^j$ is a point between $(u^j(x) + 1)$ and $(u^{j-1}(x) + 1)$. Hence let us write $c^j(x)$ as $c^j(x) = \theta^j(x) (u^j(x) + 1) + (1 - \theta^j(x)) (u^{j-1}(x) + 1)$, where $\theta^j(x) \in [0,1]$. Since $u^j(x), u^{j-1}(x) \geq 0$ and $\theta^j(x), (1 - \theta^j(x)) \in (0,1)$ we have
        \begin{equation} \vspace{-1pt}  \label{eq_lema_f''_1}
          \begin{array}{rl}
            \D {\frac{(u^j - u^{j-1})^2}{((u^j(x) + 1) + (u^{j-1}(x) + 1))^{2-s}}} & \leq \D{\frac{(u^j - u^{j-1})^2}{(c^j(x))^{2-s}}} = f''(c^j(x)) (u^j - u^{j-1})^2.
          \end{array}
        \vspace{-1pt} \end{equation} 
        
   \
        
        {\bf Step 2:} Now we estimate $\norm{(u^j + 1)^{s/2} - (u^{j-1} + 1)^{s/2}}{}^2$ by $f''(c^j(x)) (u^j - u^{j-1})^2$. If we consider the identity
        \begin{equation*} \vspace{-1pt}
          \sqrt{a} - \sqrt{b} = \frac{a - b}{\sqrt{a} + \sqrt{b}}, \ \forall a,b > 0,
        \vspace{-1pt} \end{equation*}  
        and next apply Lemma \ref{lema_T^m_elevado_a_s}, we have
        \begin{align*}
          & \norm{(u^j + 1)^{s/2} - (u^{j-1} + 1)^{s/2}}{} = \frac{\norm{(u^j + 1)^s - (u^{j-1} + 1)^s}{}}{(u^j + 1)^{s/2} + (u^{j-1} + 1)^{s/2}} \\
          & \leq s \norm{u^j - u^{j-1}}{} \frac{\big [ (u^j + 1) + (u^{j-1} + 1) \big ]^{s-1}}{(u^j + 1)^{s/2} + (u^{j-1} + 1)^{s/2}} \\
          & \leq 2^{s/2} s \norm{u^j - u^{j-1}}{} \frac{\big [ (u^j + 1) + (u^{j-1} + 1) \big ]^{s-1}}{\big [ (u^j + 1) + (u^{j-1} + 1) \big ]^{s/2}} \\
          & \leq 2^{s/2} s \norm{u^j - u^{j-1}}{} \frac{1}{\big [ (u^j + 1) + (u^{j-1} + 1) \big ]^{1 - s/2}}.
        \end{align*}
        Using \eqref{eq_lema_f''_1}
        \begin{equation} \vspace{-1pt}  \label{eq_lema_f''_2}
          \norm{(u^j + 1)^{s/2} - (u^{j-1} + 1)^{s/2}}{}^2 \leq \frac{C (u^j - u^{j-1})^2}{\big [ (u^j + 1) + (u^{j-1} + 1) \big ]^{2 - s}} \leq C f''(c^j) (u^j - u^{j-1})^2.
        \vspace{-1pt} \end{equation} 
        
  \
        
        {\bf Step 3:} Finally, we use \eqref{limitacao_diferenca_u^n_u^n-1_s<2} and \eqref{eq_lema_f''_2} to prove \eqref{convergencia_diferenca_u^n_u^n-1_s<2}. Considering the definition of $u_m^k$ and $\tilde{u}_m^k$ and using Lemma \ref{lema_T^m_elevado_a_s}, for $t \in (t^{n-1},t^n)$ we have
        \begin{align*}
          & \norm{u_m^k(t) - \tilde{u}_m^k(t)}{}^s \leq \norm{u^n - u^{n-1}}{}^s \leq \norm{((u^n + 1)^{s/2})^{2/s} - ((u^{n-1} + 1)^{s/2})^{2/s}}{}^s \\
          & \leq s \norm{(u^n + 1)^{s/2} - (u^{n-1} + 1)^{s/2}}{}^s \norm{(u^n + 1)^{s/2} + (u^{n-1} + 1)^{s/2}}{}^{2 - s}.
        \end{align*}
        Integrating and using Hölder's inequality with the conjugate powers $2/s$ and $2/(2-s)$ we obtain
        \begin{align*}
          \int_{\Omega} \norm{u_m^k(t) - \tilde{u}_m^k(t)}{}^s \ dx & \leq s \Big ( \int_{\Omega} \norm{(u^n + 1)^{s/2} - (u^{n-1} + 1)^{s/2}}{}^2 \ dx \Big )^{s/2} \\
          & \qquad \times \Big ( \int_{\Omega} \norm{(u^n + 1)^{s/2} + (u^{n-1} + 1)^{s/2}}{}^2 \ dx \Big )^{(2 - s)/2}.
        \end{align*}
        Considering the $(m,k)$-uniform bound \eqref{limitacao_u_m_L^s}, then the latter implies that
        \begin{equation*} \vspace{-1pt}
          \norma{u_m^k(t) - \tilde{u}_m^k(t)}{L^s(\Omega)}^2 \leq C \int_{\Omega} \norm{(u^n + 1)^{s/2} - (u^{n-1} + 1)^{s/2}}{}^2 \ dx.
        \vspace{-1pt} \end{equation*}  
        Then, using \eqref{eq_lema_f''_2} we obtain
        \begin{equation*} \vspace{-1pt}
          \norma{u_m^k(t) - \tilde{u}_m^k(t)}{L^s(\Omega)}^2 \leq C \int_{\Omega} f''(c^j) (u^n - u^{n-1})^2 \ dx, \mbox{ for } t \in (t_{n-1},t_n).
        \vspace{-1pt} \end{equation*}  
        Finally, if we integrate in $t$ and use \eqref{limitacao_diferenca_u^n_u^n-1_s<2} we get
        \begin{equation*} \vspace{-1pt}
          \int_0^{\infty} \norma{u_m^k(t) - \tilde{u}_m^k(t)}{L^s(\Omega)}^2 \ dt \leq C k \sum_{j=1}^{\infty} \int_{\Omega} f''(c^j) (u^j - u^{j-1})^2 \ dx \leq C k,
        \vspace{-1pt} \end{equation*}  
        which concludes the proof.
      \end{proof}
      
      Consider the relation
      \begin{equation} \vspace{-1pt}  
        \nabla u_m^k = \nabla (u_m^k + 1) = \nabla \big ( (u_m^k + 1)^{s/2} \big )^{2/s} = \dfrac{2}{s} (u_m^k + 1)^{1 - s/2} \ \nabla (u_m^k + 1)^{s/2}.
        \label{gradiente_u_m_em_funcao_de_gradiente_u_m^s/2}
      \vspace{-1pt} \end{equation} 
      Taking into account that we are considering $s \in [1,2)$, we can use \eqref{limitacao_u_m_L^s} to obtain 
      \begin{equation*} \vspace{-1pt}
        (u_m^k + 1)^{1-s/2} \mbox{ is bounded in } L^{\infty}(0,\infty;L^{2s/(2-s)}(\Omega))
      \vspace{-1pt} \end{equation*}  
      and then \eqref{limitacao_nabla_u_m_elevado_a_s_sobre_2} and \eqref{gradiente_u_m_em_funcao_de_gradiente_u_m^s/2} to conclude that 
      \begin{equation} \vspace{-1pt}  
        \nabla u_m^k \mbox{ is bounded in } L^2(0,\infty;L^s(\Omega)).
        \label{limitacao_nabla_u_m_L^s}
      \vspace{-1pt} \end{equation} 
      In conclusion, using \eqref{limitacao_u_m_L^s_2}, \eqref{limitacao_nabla_u_m_L^s} and the Poincare's inequality for zero mean functions (Lemma \ref{lema_desigualdade_poincare_media_nula}),
      \begin{equation} \vspace{-1pt}  
        u_m^k - u^\ast \mbox{ is bounded in } L^{\infty}(0,\infty;L^s(\Omega)) \cap L^2(0,\infty;W^{1,s}(\Omega)).
        \label{limitacao_W_1_s}
      \vspace{-1pt} \end{equation} 
      
      Considering the chemotaxis term, we write $T^m(u_m^k) \nabla (z_m^k)^2$ as
  \begin{equation*} \vspace{-1pt}
    T^m(u_m^k) \nabla (z_m^k)^2 = 2 T^m(u_m^k)^{1-s/2} T^m(u_m^k)^{s/2} z_m^k \nabla z_m^k.
  \vspace{-1pt} \end{equation*}  
  Then, because of \eqref{limitacao_u_m_L^s}, we have $T^m(u_m^k)^{1-s/2}$ bounded in $L^{\infty}(0,\infty;L^{2s/(2-s)}(\Omega))$ and , because of \eqref{limitacao_aux_Delta_z_m_s_intermediario}, $T^m(u_m^k)^{s/2} z_m^k \nabla z_m^k$ is bounded in $L^2(0,\infty;L^2(\Omega))$. Hence we conclude that
  \begin{equation} \vspace{-1pt}  
    T^m(u_m^k) \nabla (z_m^k)^2 \mbox{ is bounded in } L^2(0,\infty;L^s(\Omega)).
    \label{limitacao_termo_chemotaxis_s}
  \vspace{-1pt} \end{equation} 
  Then, if we consider the $u$-equation of \eqref{NLTD_equiv}, from \eqref{limitacao_W_1_s} and \eqref{limitacao_termo_chemotaxis_s} we have 
  \begin{equation*} \vspace{-1pt}
    \partial_t \tilde{u}_m^k \mbox{ is bounded in } L^2 \Big (0,\infty;\big ( W^{1,s/(s-1)}(\Omega) \big )' \Big ).
  \vspace{-1pt} \end{equation*}  
      
      Now we turn to the $z$-equation of \eqref{NLTD_equiv}, rewritten as
      \begin{equation} \vspace{-1pt}  
        \begin{array}{l}
          \partial_t \tilde{z}_m^k - \dfrac{\norm{\nabla z_m^k}{}^2}{z_m^k} - \Delta z_m^k + T^m(u^{\ast})^s (z_m^k - \dfrac{\alpha^2}{z_m^k}) \\
          = - \dfrac{1}{2} (T^m(u_m^k)^s - T^m(u^{\ast})^s) (z_m^k - \dfrac{\alpha^2}{z_m^k}).
        \end{array}
        \label{equacao_z_m_equivalente_s<2}
      \vspace{-1pt} \end{equation} 
      Since $z_m^k \geq \alpha$, we have $z_m^k \geq z_m^k - \alpha^2/z_m^k \geq z_m^k - \alpha \geq 0$ and then we can write
      \begin{equation} \vspace{-1pt}  
        \begin{array}{l}
          \partial_t \tilde{z}_m^k - \dfrac{\norm{\nabla (z_m^k - \alpha)}{}^2}{z_m^k} - \Delta (z_m^k - \alpha) + T^m(u^{\ast})^s (z_m^k - \alpha) \\
          \leq \dfrac{1}{2} \norm{T^m(u_m^k)^s - T^m(u^{\ast})^s)}{} z_m^k.
        \end{array}
        \label{eq_z_m_equivalente_s<2}
      \vspace{-1pt} \end{equation} 
      Analyzing the term on the right hand side of \eqref{eq_z_m_equivalente_s<2}, we have
      \begin{equation} \vspace{-1pt}  
        T^m(u_m^k)^s - T^m(u^{\ast})^s \mbox{ is bounded in } L^2(0,\infty;L^{3/2}(\Omega)).
        \label{limitacao_u_s_menos_u_estrela_s}
      \vspace{-1pt} \end{equation} 
      In fact, using Lemma \ref{lema_T^m_elevado_a_s}, we get
      \begin{align*}
        \norm{T^m(u_m^k)^s - T^m(u^{\ast})^s}{} \leq s^{3/2} \norm{T^m(u_m^k) + T^m(u^{\ast})}{}^{s-1} \norm{u_m^k - u^{\ast}}{}
      \end{align*}
      and therefore
      \begin{align*}
        & \int_{\Omega} \norm{T^m(u_m^k)^s - T^m(u^{\ast})^s}{}^{3/2} \ dx \leq s^{3/2} \int_{\Omega} \norm{T^m(u_m^k) + T^m(u^{\ast})}{}^{3(s-1)/2} \ dx \int_{\Omega} \norm{u_m^k - u^{\ast}}{}^{3/2} \ dx.
      \end{align*}
      Using Hölder's inequality with the conjugate exponents $2s/(3s-3)$ and $2s/(3-s)$, we obtain
      \begin{align*}
        & \norma{T^m(u_m^k)^s - T^m(u^{\ast})^s}{L^{3/2}(\Omega)}^2 \leq s^2 \norma{T^m(u_m^k) + T^m(u^{\ast})}{L^s(\Omega)}^{2(s-1)} \norma{u_m^k - u^{\ast}}{L^{3s/(3-s)}(\Omega)}^2.
      \end{align*}
      Then, recalling the $(m,k)$-uniform bound \eqref{limitacao_W_1_s} and the Sobolev embedding $W^{1,s}(\Omega) \subset L^{3s/(3-s)}(\Omega)$ we obtain \eqref{limitacao_u_s_menos_u_estrela_s}.
      
      With this information, now we can test \eqref{eq_z_m_equivalente_s<2} by $k \ (z_m^k - \alpha)$ and get, for each time interval $(t_{n-1},t_n)$,
      \begin{align*}
        & \frac{1}{2} ( \norma{z^n - \alpha}{L^2(\Omega)}^2 - \norma{z^{n-1} - \alpha}{L^2(\Omega)}^2 ) + k \norma{\nabla (z^n - \alpha)}{L^2(\Omega)}^2 + k \frac{T^m(u^{\ast})^s}{2} \norma{z^n - \alpha}{L^2(\Omega)}^2 \\
        & \leq C k \int_{\Omega}{\norm{T^m(u^n)^s - T^m(u^{\ast})^s}{} z^n (z^n - \alpha) \ dx} + \int_{\Omega} \frac{\norm{\nabla z^n}{}^2}{z^n} (z^n - \alpha) \ dx \\
        & \leq C k \norma{z^0}{L^{\infty}(\Omega)} \norma{T^m(u^n)^s - T^m(u^{\ast})^s}{L^{3/2}(\Omega)} \norma{z^n - \alpha}{L^3(\Omega)} + k \norma{\nabla z^n}{L^2(\Omega)}^2 \\
        & \leq C(\delta) k \norma{T^m(u^n)^s - T^m(u^{\ast})^s}{L^{3/2}(\Omega)}^2 + \delta k \norma{z^n - \alpha}{L^2(\Omega)}^2 + (1 + \delta) k \norma{\nabla z^n}{L^2(\Omega)}^2.
      \end{align*}
      Note that if $u_0 \not\equiv 0$ then have $T^m(u^{\ast}) = u^{\ast} > 0$, for all $m \geq u^{\ast}$. Hence, choosing $\delta > 0$ small enough, we conclude that, for $m \geq u^{\ast}$, there is $\beta > 0$ such that
      \begin{align*}
        & \frac{1}{2} ( \norma{z^n - \alpha}{L^2(\Omega)}^2 - \norma{z^{n-1} - \alpha}{L^2(\Omega)}^2 ) + k \beta \norma{z^n - \alpha}{L^2(\Omega)}^2 \\
        & \leq C k \norma{T^m(u_m)^s - T^m(u^{\ast})^s}{L^{3/2}(\Omega)}^2 + k \norma{\nabla z^n}{L^2(\Omega)}^2.
      \end{align*}
      Therefore, summing in $n$ and using \eqref{limitacao_u_s_menos_u_estrela_s} and Lemma \ref{lema_limitacao_grad_z_uniforme_wrt_k_m_T} we obtain
      \begin{equation} \vspace{-1pt}  
        z_m^k - \alpha \mbox{ is bounded in } L^2(0,\infty;L^2(\Omega)).
        \label{estimativa_z_m_L^2}
      \vspace{-1pt} \end{equation} 
      Hence, in view of \eqref{limitacao_aux_Dz_m_s_intermediario}, \eqref{limitacao_aux_Delta_z_m_s_intermediario} and \eqref{estimativa_z_m_L^2} we have
      \begin{equation*} \vspace{-1pt}
        z_m^k - \alpha \mbox{ is bounded in } L^2(0,\infty;H^2(\Omega)).
      \vspace{-1pt} \end{equation*}  
      
      With the $m$-uniform bounds obtained so far we can derive a $m$-uniform bound for $\partial_t \tilde{z}_m^k$ in $L^2(0,\infty;L^{3/2}(\Omega))$. In fact, we notice that
      \begin{equation*} \vspace{-1pt}
        0 \leq z_m^k - \frac{\alpha^2}{z_m^k} = \frac{(z_m^k)^2 - \alpha^2}{z_m^k} \leq \frac{(z_m^k)^2 - \alpha^2}{\alpha} \leq \frac{\norma{z^0}{L^{\infty}(\Omega)} + \alpha^2}{\alpha} (z_m^k - \alpha).
      \vspace{-1pt} \end{equation*}  
      Hence, accounting for \eqref{estimativa_z_m_L^2}, we also have
      \begin{equation} \vspace{-1pt}  \label{estimativa_z_m_L^2_consequencia}
        z_m^k - \frac{\alpha^2}{z_m^k} \mbox{ is bounded in } L^2(0,\infty;L^2(\Omega)).
      \vspace{-1pt} \end{equation} 
      Therefore, going back to \eqref{equacao_z_m_equivalente_s<2}, recalling that $z_m^k$ is uniformly bounded in $L^{\infty}(0,\infty;L^{\infty}(\Omega))$ with respect to $m$ and $k$ and considering \eqref{estimativa_z_m_L^2_consequencia} and \eqref{limitacao_u_s_menos_u_estrela_s}, we conclude that
      \begin{equation*} \vspace{-1pt}
        \partial_t \tilde{z}_m^k \mbox{ is bounded in } L^2(0,\infty;L^{3/2}(\Omega)).
      \vspace{-1pt} \end{equation*}  
      
      Now we are going to obtain compactness for $\{ \tilde{u}_m^k \}$ and $\{ u_m^k \}$, which is necessary in order to pass to the limit as $m \to \infty$ and $k \to 0$ in the nonlinear terms of \eqref{NLTD}. Because of \eqref{limitacao_u_m_L^s} and \eqref{limitacao_nabla_u_m_elevado_a_s_sobre_2}, we have that
      \begin{equation*} \vspace{-1pt}
        (u_m^k + 1)^{s/2} \mbox{ is bounded in } L^{\infty}(0,\infty;L^2(\Omega)) \cap L^2(0,T;H^1(\Omega)),
      \vspace{-1pt} \end{equation*}  
      for every finite $T > 0$. Using the Sobolev inequality $H^1(\Omega) \subset L^6(\Omega)$ and interpolation inequalities we obtain
      \begin{equation*} \vspace{-1pt}
        (u_m^k)^{s/2} \mbox{ is bounded in } L^{10/3}(0,T;L^{10/3}(\Omega)),
      \vspace{-1pt} \end{equation*}  
      which is equivalent to
      \begin{equation} \vspace{-1pt}  \label{limitacao_u_m_L_5s/3}
        u_m^k \mbox{ is bounded in } L^{5s/3}(0,T;L^{5s/3}(\Omega)).
      \vspace{-1pt} \end{equation} 
      By using \eqref{limitacao_u_m_L^s} and \eqref{limitacao_u_m_L_5s/3} in  \eqref{gradiente_u_m_em_funcao_de_gradiente_u_m^s/2} (remind that $s \in [1,2)$), we also have
      \begin{equation} \vspace{-1pt}  
        u_m^k \mbox{ is bounded in } L^{5s/(3+s)}(0,T;W^{1,5s/(3+s)}(\Omega)).
       \label{limitacao_W_5s_3+s}
      \vspace{-1pt} \end{equation} 
      For any norm $\norma{\cdot}{}$ we have
      \begin{equation*} \vspace{-1pt}
        \norma{\tilde{u}_m^k(t)}{} \leq \norma{\tilde{u}_m^k(t) - u_m^k(t)}{} + \norma{u_m^k(t)}{},
      \vspace{-1pt} \end{equation*}  
      \begin{equation*} \vspace{-1pt}
        \norma{\tilde{u}_m^k(t) - u_m^k(t)}{} \leq \norma{u^n - u^{n-1}}{} \leq \norma{u^n}{} + \norma{u^{n-1}}{}, \ \forall t \in (t_{n-1}, t_n).
      \vspace{-1pt} \end{equation*}  
      Thus, because of \eqref{limitacao_u_m_L_5s/3} and \eqref{limitacao_W_5s_3+s} we also conclude that
      \begin{equation} \vspace{-1pt}  
        \tilde{u}_m^k \mbox{ is bounded in } L^{5s/3}(0,T;L^{5s/3}(\Omega)) \cap L^{5s/(3+s)}(0,T;W^{1,5s/(3+s)}(\Omega)).
      \vspace{-1pt} \end{equation} 
      
      We observe that $W^{1,5s/(3+s)}(\Omega) \subset L^q(\Omega)$, with continuous embedding for $q = 15s/(9-2s)$ and compact embedding for $q \in [1,15s/(9-2s))$. Then, since $s \in [1,2)$, we have $5s/3 < 15s/(9-2s)$ and therefore the embedding $W^{1,5s/(3+s)}(\Omega) \subset L^{5s/3}(\Omega)$ is compact. Note also that $q = 5s/3 \geq 5/3 > 1$. Applying Lemma \ref{lema_Simon} with 
      \begin{equation*} \vspace{-1pt}
        X = W^{1,5s/(3+s)}(\Omega), \hspace{1cm} B = L^{5s/3}(\Omega), \hspace{1cm} Y = \big ( H^3(\Omega) \big )',
      \vspace{-1pt} \end{equation*}  
      we conclude that there is a subsequence of $\{ \tilde{u}_m^k \}$ (still denoted by $\{ u_m \}$) and a limit function $u$ such that
      \begin{equation*} \vspace{-1pt}
        \begin{array}{c}
          \tilde{u}_m^k \longrightarrow u \mbox{ weakly in } L^{5s/(3+s)}(0,T;W^{1,5s/(3+s)}(\Omega)), \ \forall T > 0, \\
          \tilde{u}_m^k \longrightarrow u \mbox{ strongly in } L^p(0,T;L^p(\Omega)), \ \forall p \in [1,5s/3), \ \forall T > 0.
        \end{array}
      \vspace{-1pt} \end{equation*}  
      We note that, because of the $(m,k)$-independent bounds obtained for $u_m^k$ and the convergence \eqref{convergencia_diferenca_u^n_u^n-1_s<2}, these convergences are also valid if we replace $\tilde{u}_m^k$ by $u_m^k$. Therefore we have
      \begin{equation*} \vspace{-1pt}
        u_m^k \longrightarrow u \mbox{ weakly in } L^{5s/(3+s)}(0,T;W^{1,5s/(3+s)}(\Omega)), \ \forall T > 0,
      \vspace{-1pt} \end{equation*}  
      \begin{equation} \vspace{-1pt}  \label{convergencia_u_m_s_intermediario}
        u_m^k \longrightarrow u \mbox{ strongly in } L^p(0,T;L^p(\Omega)), \ \forall p \in [1,5s/3), \ \forall T > 0.
      \vspace{-1pt} \end{equation} 
      
      Using the Dominated Convergence Theorem we conclude from \eqref{convergencia_u_m_s_intermediario} that
      \begin{equation} \vspace{-1pt}  \label{convergencia_T^m_u_s_intermediario}
        T^m(u_m^k) \rightarrow u \mbox{ strongly in } L^p(0,T;L^p(\Omega)), \ \forall p \in [1,5s/3), \ \forall T > 0.
      \vspace{-1pt} \end{equation} 
      It stems from the convergence \eqref{convergencia_T^m_u_s_intermediario} and Lemma \ref{lema_convergencia_w_elevado_a_s} that
      \begin{equation} \vspace{-1pt}  \label{convergencia_T^m_u_elevado_a_s_s_intermediario}
        (T^m(u_m^k))^s \rightarrow u^s \mbox{ strongly in } L^q(0,T;L^q(\Omega)), \ \forall q \in [1,5/3), \ \forall T > 0.
      \vspace{-1pt} \end{equation} 
      
      The compactness of $\{ z_m^k \}$ is also necessary. Concerning the functions $z_m^k$ and $\tilde{z}_m^k$, if we consider the $(m,k)$-independent bounds derived so far, \eqref{convergencia_diferenca_z_m_em_H_1} and use the compactness result Lemma \ref{lema_Simon}, then we conclude that there are subsequences of $\{ z_m^k \}$ and $\{ \tilde{z}_m^k \}$ (still denoted by $\{ z_m^k \}$ and $\{ \tilde{z}_m^k \}$) and a limit function $z$ such that, for each $T > 0$,
      \begin{equation} \vspace{-1pt}  \label{convergencia_z_m_s_intermediario}
        \begin{array}{c}
          z_m^k \rightarrow z \mbox{ weakly* in } L^{\infty}(0,\infty;L^{\infty}(\Omega)) \cap L^{\infty}(0,\infty;H^1(\Omega)), \\
          z_m^k \rightarrow z \mbox{ weakly in } L^2(0,\infty;H^2(\Omega)), \\
          z_m^k \rightarrow z \mbox{ strongly in } L^2(0,T;H^1(\Omega)) \cap L^p(0,T;L^p(\Omega)), p \in [1, \infty), \\
          \nabla z_m^k \rightarrow \nabla z \mbox{ weakly in } L^4(0,\infty;L^4(\Omega)), \\
          \mbox{and } \partial_t \tilde{z}_m^k \rightarrow \partial_t z \mbox{ weakly in } L^2(0,\infty;L^{3/2}(\Omega)).
        \end{array}
      \vspace{-1pt} \end{equation} 
      Taking the nonlinear terms of \eqref{NLTD_equiv}, where we have to pass to the limit as $(m,k) \to (\infty,0)$, it is convenient to consider two functions of $z_m^k$, namely $g_1(z_m^k) = (z_m^k)^2$ and $g_2(z_m^k) = 1/z_m^k$. Since $0 < \alpha \leq z_m^k \leq \norma{z^0}{L^{\infty}(\Omega)}$, we are able to use the $(m,k)$-independent bounds derived so far and \eqref{convergencia_z_m_s_intermediario} to show that
      \begin{equation} \vspace{-1pt}  
        \begin{array}{c}
          1/z_m^k \rightarrow 1/z \mbox{ strongly in } L^p(0,T;L^p(\Omega)), \mbox{ for each } T > 0 \mbox{ and } p \in [1, \infty), \\
          \norm{\nabla z_m^k}{}^2 \rightarrow \norm{\nabla z}{}^2 \mbox{ weakly in } L^2(0,\infty;L^2(\Omega)), \\
          \nabla (z_m^k)^2 \rightarrow \nabla (z)^2 \mbox{ weakly in } L^4(0,\infty;L^4(\Omega)).
        \end{array}
        \label{convergencia_z_m_s_intermediario_2}
      \vspace{-1pt} \end{equation} 
      
      Now we are going to use the weak and strong convergences that we proved for $u_m^k$, $\tilde{u}_m^k$, $z_m^k$ and $\tilde{z}_m^k$ to pass to the limit as $(m,k) \to (\infty,0)$ in the equations of problem \eqref{NLTD}. We are going to identify the limits of the nonlinear terms
      \begin{equation*}
        T^m(u_m^k) \nabla (z_m^k)^2, \ \dfrac{\norm{\nabla z_m^k}{}^2}{z_m^k}, \ T^m(u_m^k)^s z_m^k \ \mbox{ and } \ \dfrac{T^m(u_m^k)^s}{z_m^k},
      \end{equation*}
      with $u \nabla (z)^2$, $\dfrac{\norm{\nabla z}{}^2}{z}$, $u^s z$ and $\dfrac{u^s}{z}$, respectively. In fact, considering the chemotaxis term, because of \eqref{convergencia_T^m_u_s_intermediario}, \eqref{limitacao_aux_Dz_m_s_intermediario} and \eqref{convergencia_z_m_s_intermediario_2}, we have
      \begin{equation*} \vspace{-1pt}
        T^m(u_m^k) \nabla (z_m^k)^2 \longrightarrow u \nabla (z)^2 \mbox{ weakly in } L^{20s/(5s+12)}(0,T;L^{20s/(5s+12)}(\Omega)),
      \vspace{-1pt} \end{equation*}  
      for each $T > 0$. Using \eqref{convergencia_z_m_s_intermediario} and \eqref{convergencia_z_m_s_intermediario_2} we also conclude that
      \begin{equation*} \vspace{-1pt}
        \frac{\norm{\nabla z_m^k}{}^2}{z_m^k} \longrightarrow \frac{\norm{\nabla z}{}^2}{z} \mbox{ weakly in } L^2(0,T;L^2(\Omega)), \mbox{ for each } T > 0.
      \vspace{-1pt} \end{equation*}  
      Regarding the consumption term, considering \eqref{convergencia_T^m_u_elevado_a_s_s_intermediario} and \eqref{convergencia_z_m_s_intermediario}, we prove that
      \begin{equation*} \vspace{-1pt}
        T^m(u_m^k)^s z_m^k \longrightarrow u^s z \mbox{ weakly in } L^{5/3}(0,T;L^{5/3}(\Omega)), \mbox{ for each } T > 0.
      \vspace{-1pt} \end{equation*}  
      Finally, using \eqref{convergencia_T^m_u_elevado_a_s_s_intermediario} and \eqref{convergencia_z_m_s_intermediario_2}, we obtain
      \begin{equation*} \vspace{-1pt}
        \frac{T^m(u_m^k)^s}{z_m^k} \longrightarrow \frac{u^s}{z} \mbox{ weakly in } L^{5/3}(0,T;L^{5/3}(\Omega)), \mbox{ for each } T > 0.
      \vspace{-1pt} \end{equation*}  
  
      With these identifications and all previous convergences, it is possible to pass to the limit as $(m,k) \to (\infty,0)$ in each term of the equations of \eqref{NLTD}. This concludes the proof of existence of a solution $(u,z)$ of \eqref{problema_P_u_z} for $s \in [1,2)$. In order to finish, we will prove the regularity (up to infinite time) of $u$.
  
      From \eqref{limitacao_u_m_L^s} and \eqref{limitacao_nabla_u_m_elevado_a_s_sobre_2} there exists a subsequence of $\{ (\tilde{u}_m^k + 1)^{s/2} \}$, still denoted by $\{ (\tilde{u}_m^k + 1)^{s/2} \}$, and a limit function $\varphi$ such that
      \begin{equation*} \vspace{-1pt}
        \begin{array}{c}
          (u_m^k + 1)^{s/2} \longrightarrow \varphi \mbox{ \ weakly* in } L^{\infty}(0,\infty;L^2(\Omega)) \\
          \nabla (u_m^k + 1)^{s/2} \longrightarrow \nabla \varphi \mbox{ \ weakly in } L^2(0,\infty;L^2(\Omega)).
        \end{array}
      \vspace{-1pt} \end{equation*}  
      Then, using the strong convergence \eqref{convergencia_u_m_s_intermediario}, the continuity of the function $u_m^k \mapsto f(u_m^k) = (u_m^k + 1)^{s/2}$ and the Dominated Convergence Theorem, we prove that $\varphi = (u + 1)^{s/2}$.
  
      Analogously, because of \eqref{limitacao_aux_Delta_z_m_s_intermediario}, up to a subsequence, there is a limit function $\phi$ such that
      \begin{equation*} \vspace{-1pt}
        T^m(u_m^k)^{s/2} \nabla z_m^k \longrightarrow \phi \mbox{ weakly in } L^2(0,\infty;L^2(\Omega)).
      \vspace{-1pt} \end{equation*}  
      Using the convergences \eqref{convergencia_T^m_u_elevado_a_s_s_intermediario} and \eqref{convergencia_z_m_s_intermediario}, we prove that $\phi = u^{s/2} \nabla z$.
  
      Therefore we have proved the global in time regularity
      \begin{equation} \vspace{-1pt}  
        \begin{array}{c}
          (u + 1)^{s/2} \in L^{\infty}(0,\infty;L^2(\Omega)), \quad
          \nabla (u + 1)^{s/2} \in L^2(0,\infty;L^2(\Omega)), \\
          u^{s/2} \nabla z \in L^2(0,\infty;L^2(\Omega)).
        \end{array}
        \label{regularidade_termos_limite}
      \vspace{-1pt} \end{equation} 
      Considering \eqref{regularidade_termos_limite} and proceeding as in the obtaining of \eqref{limitacao_W_1_s} and \eqref{limitacao_termo_chemotaxis_s} we conclude the global in time regularity $u \in L^{\infty}(0,\infty;L^s(\Omega))$, $\nabla u, \ u \nabla (z)^2 \in L^2(0,\infty;L^s(\Omega))$. This finishes the proof that $(u,z)$ is a weak solution of \eqref{problema_P_u_z} and that $\{ (u_m^k,z_m^k) \}$ converges to $(u,z)$ as $(m,k) \to (\infty, 0)$ in the sense indicated in this section, for $s \in [1,2)$.

      
      \subsection{Estimates for \texorpdfstring{$\boldsymbol{(u_m^k,z_m^k)}$}{(u\_m\^\{k,r\},z\_m\^\{k,r\})} and passage to the limit for \texorpdfstring{$\boldsymbol{s \geq 2}$}{s geq 2}}
      
      \label{subsec:estimativas_u_m_s_geq_2}
      
        The procedure for the case $s \geq 2$ is much more similar to the case $s \geq 2$ in \cite{ViannaGuillen2023uniform}. In the sequel, we highlight the main steps of the proof and refer the reader to \cite{ViannaGuillen2023uniform} for details. First we note that, multiplying the energy inequality \eqref{estimativa_u_v_m_1_s_geq_2} from Lemma \ref{lemma_estimativa_u_v_m_1_s_geq_2} by $k$ and summing in $n$, for each $n \in \mathbb{N}$, we have
        \begin{equation} \vspace{-1pt}  \label{limitacao_nabla_T^m_u_m_s_geq_2}
          \nabla T^m(u_m^k)^{s/2} \mbox{ is bounded in } L^2(0,\infty;L^2(\Omega)),
        \vspace{-1pt} \end{equation} 
        \begin{equation} \vspace{-1pt}  \label{limitacao_T^m_u_m_s_geq_2}
          T^m(u_m^k)^{s/2} \mbox{ is bounded in } L^{\infty}(0,\infty;L^2(\Omega)).
        \vspace{-1pt} \end{equation} 
        From \eqref{limitacao_T^m_u_m_s_geq_2} and \eqref{limitacao_nabla_T^m_u_m_s_geq_2} we conclude that
        \begin{equation} \vspace{-1pt}  \label{limitacao_T^m_u_m_5s/3_s_geq_2}
          T^m(u_m^k) \mbox{ is bounded in } L^{5s/3}(0,T;L^{5s/3}(\Omega)).
        \vspace{-1pt} \end{equation} 
        
        Analogously to \cite{ViannaGuillen2023uniform}, we use \eqref{limitacao_aux_Delta_z_m_s_intermediario} and Lemma \ref{lema_limitacao_grad_z_uniforme_wrt_k_m_T} to prove that
        \begin{equation} \vspace{-1pt}  
          T^m(u_m^k) \nabla (z_m^k)^2 \mbox{ is bounded in } L^2(0,\infty;L^2(\Omega)).
          \label{limitacao_termo_misto_u_v_m_s_geq_2}
        \vspace{-1pt} \end{equation} 
      
        Now we can test the $u^n$-equation of \eqref{NLTD} by $k u^n$ and, after bounding some terms, we sum the resulting inequality in $n \in \mathbb{N}$ and use \eqref{limitacao_termo_misto_u_v_m_s_geq_2} to conclude that
        \begin{equation} \vspace{-1pt}  
          u_m^k \mbox{ is bounded in } L^{\infty}(0,\infty;L^2(\Omega)),
          \label{limitacao_u_m_s_geq_2}
        \vspace{-1pt} \end{equation} 
        \begin{equation} \vspace{-1pt}  
          \nabla u_m^k \mbox{ is bounded in } L^2(0,\infty;L^2(\Omega)),
          \label{limitacao_nabla_u_m_s_geq_2}
        \vspace{-1pt} \end{equation} 
        \begin{equation*} \vspace{-1pt}
          \D{\sum_{j=1}^{\infty}}{\norma{u^j - u^{j-1}}{L^2(\Omega)}^2} \leq C.
          \label{limitacao_diferenca_u_L_2}
        \vspace{-1pt} \end{equation*}  
        Analogously to \eqref{convergencia_diferenca_z_m_em_H_1}, using \eqref{limitacao_diferenca_u_L_2}, we prove that
        \begin{equation} \vspace{-1pt}  
          \norma{u_m^k - \tilde{u}_m^k}{L^2(0,\infty;L^2(\Omega))}^2 \leq C k.
          \label{convergencia_diferenca_u_m_em_L_2}
        \vspace{-1pt} \end{equation} 
        
        Then, if we consider the $u$-equation of \eqref{NLTD_equiv}, by applying \eqref{limitacao_nabla_u_m_s_geq_2} and \eqref{limitacao_termo_misto_u_v_m_s_geq_2} we conclude that
        \begin{equation} \vspace{-1pt}  
          \partial_t \tilde{u}_m^k \mbox{ is bounded in } L^2(0,\infty;(H^1(\Omega))').
          \label{limitacao_u_m_derivada_tempo_s_geq_2}
        \vspace{-1pt} \end{equation} 
        
        Using the $(m,k)$-independent bounds obtained so far, we can follow the ideas of \cite{ViannaGuillen2023uniform} and subsection \ref{subsec:estimativas_u_m_s<2} in order to prove that
        \begin{equation} \vspace{-1pt}  
          z_m^k - \alpha \mbox{ is bounded in } L^2(0,\infty;L^2(\Omega)),
          \label{estimativa_z_m_L^2_s_geq_2}
        \vspace{-1pt} \end{equation} 
        \begin{equation*} \vspace{-1pt}
          \partial_t \tilde{z}_m^k \mbox{ is bounded in } L^2(0,\infty;L^{3/2}(\Omega)).
        \vspace{-1pt} \end{equation*}  
      
        Now, using \eqref{limitacao_u_m_s_geq_2}, \eqref{limitacao_nabla_u_m_s_geq_2}, \eqref{convergencia_diferenca_u_m_em_L_2} and \eqref{limitacao_u_m_derivada_tempo_s_geq_2} we conclude that there are subsequences of $\{ u_m^k \}$ and $\{ \tilde{u}_m^k \}$, still denoted by $\{ u_m^k \}$ and $\{ \tilde{u}_m^k \}$, and a limit function $u$ such that
        \begin{equation*} \vspace{-1pt}
          \begin{array}{c}
            u_m^k \longrightarrow u \mbox{ weakly* in } L^{\infty}(0,\infty;L^2(\Omega)), \\
            \nabla u_m^k \longrightarrow \nabla u \mbox{ weakly in } L^2(0,\infty;L^2(\Omega)), \\
            \partial_t \tilde{u}_m^k \longrightarrow u \mbox{ weakly in } L^2 \Big ( 0,\infty;\big (H^1(\Omega) \big )' \Big ).
          \end{array}
        \vspace{-1pt} \end{equation*}  
        By applying the compactness result Lemma \ref{lema_Simon} to the sequence $\{ \tilde{u}_m^k \}$ and using \eqref{convergencia_diferenca_u_m_em_L_2} we have
        \begin{equation*} \vspace{-1pt}
          u_m^k \longrightarrow u \mbox{ strongly in } L^2(0,T;L^2(\Omega)), \ \forall T > 0.
        \vspace{-1pt} \end{equation*}  
        Using the Dominated Convergence Theorem and \eqref{limitacao_T^m_u_m_5s/3_s_geq_2} we can also prove that
        \begin{equation*} \vspace{-1pt}
          T^m(u_m^k) \longrightarrow u \mbox{ strongly in } L^p(0,T;L^p(\Omega)), \ \forall p \in (1,5s/3),
        \vspace{-1pt} \end{equation*}  
        and then Lemma \ref{lema_convergencia_w_elevado_a_s} yields
        \begin{equation*} \vspace{-1pt}
          T^m(u_m^k)^s \longrightarrow u^s \mbox{ strongly in } L^p(0,T;L^p(\Omega)), \ \forall p \in (1,5/3).
        \vspace{-1pt} \end{equation*}  
        
        From the global in time estimate \eqref{limitacao_T^m_u_m_s_geq_2} we conclude that, up to a subsequence,
          \begin{equation*} \vspace{-1pt}
            T^m(u_m^k) \to u \mbox{ weakly* in } L^{\infty}(0,\infty;L^s(\Omega)),
         \vspace{-1pt} \end{equation*}  
         hence, in particular, $u \in L^{\infty}(0,\infty;L^s(\Omega))$.
         
         For $s \geq 2$, if we consider the functions $z_m^k$ and $\tilde{z}_m^k$, we have the same $(m,k)$-independent estimates that we had for $s \in [1,2)$. Then we have the same convergences given in \eqref{convergencia_z_m_s_intermediario} and \eqref{convergencia_z_m_s_intermediario_2}.
         
         Following the ideas of Subsection \ref{subsec:estimativas_u_m_s<2}, we are able identify the limits of 
         \begin{equation*}
           T^m(u_m^k) \nabla (z_m^k)^2, \ \norm{\nabla z_m^k}{}^2/z_m^k, \ T^m(u_m^k)^s z_m^k \ \mbox{ and } \ T^m(u_m^k)^s/z_m^k
         \end{equation*}
         with $u \nabla (z)^2$, $\norm{\nabla z}{}^2/z$, $u^s z$ and $u^s/z$, respectively.
         
         This finishes the proof that $(u,z)$ is a solution of \eqref{problema_P_u_z} and that $\{ (u_m^k,z_m^k) \}$ converges to $(u,z)$ as $(m,k) \to (\infty, 0)$ in the sense indicated in this section, for $s \geq 2$.

        
        \section{Convergence of \texorpdfstring{$\boldsymbol{v_m^k}$}{v\_m\^\{k,r\}}}
          \label{subsec:convergencia_de_v_m^kr}
        
        Until this point, for any $s \geq 1$ fixed, we have proved that $(u_m^k, z_m^k)$ converges to a weak solution $(u,z)$ of \eqref{problema_P_u_z} as $(m,k) \to (\infty, 0)$. Now, to conclude the proof of Theorem \ref{teo_principal}, we are going to prove that $(u_m^k, v_m^k)$ converges to a weak solution $(u,v)$ of \eqref{problema_P} as $(m,k) \to (\infty, 0)$, where $v_m^k$ is given  either by  \eqref{v-z} or by \eqref{v-u}. For simplicity, we consider $(u_m^k,z_m^k)$ to be the subsequence which converges to the limit function $(u,z)$.
        
        
        \subsection{\texorpdfstring{$\boldsymbol{v_m^k}$}{v\_m\^\{k,r\}} given by \texorpdfstring{\eqref{v-z}}{{\bf v\_m\^\{k,r\} = (z\_m\^\{k,r\})\^2 - alpha\^2}}}

          In this case, it suffices to show that the sequence $v_m^k = (z_m^k)^2 - \alpha^2$ converges to $v = z^2 - \alpha^2$ as $(m,k) \to (\infty, 0)$. Then, thanks to the equivalence of problems \eqref{problema_P_u_z} and \eqref{problema_P} (see Lemma \ref{lema_equivalencia_problemas}), we know that $(u,z)$ is a solution of \eqref{problema_P_u_z} if, and only if, $(u,v)$ is a solution of \eqref{problema_P}, with $v = z^2 - \alpha^2$.
          
          In fact, if we consider the $(m,k)$-uniform bounds obtained for $z_m^k$ (especially the pointwise estimates of Theorem \ref{teo_existencia_NLTD}, \eqref{limitacao_aux_Dz_m_s_intermediario} and \eqref{limitacao_aux_Delta_z_m_s_intermediario}) and the convergences listed in \eqref{convergencia_z_m_s_intermediario}, we can prove by straightforward calculations that the sequence $v_m^k = (z_m^k)^2 - \alpha^2$ converges to $v = z^2 - \alpha^2$ in the same senses indicated in \eqref{convergencia_z_m_s_intermediario}. Hence, by Lemma \ref{lema_equivalencia_problemas}, we have that $(u,v)$ is a solution of \eqref{problema_P}.
       
          Therefore, we conclude that $(u_m^k, v_m^k)$ converges to $(u,v)$, a solution of \eqref{problema_P}, as $(m,k) \to (\infty, 0)$.

        
        \subsection{\texorpdfstring{$\boldsymbol{v_m^k}$}{v\_m\^\{k,r\}} given \texorpdfstring{by \eqref{v-u}}{{\bf in terms of u\_m\^\{k,r\}}}}
          
          We rewrite \eqref{v-u} as
          \begin{equation} \vspace{-1pt}  
            \partial_t \tilde{v}_m^k - \Delta v_m^k + T^m(u_m^k)^s v_m^k = 0.
            \label{eq_de_v_m^kr}
          \vspace{-1pt} \end{equation} 
          Taking the $(m,k)$-independent bounds obtained for $T^m(u_m^k)$ in Subsections \ref{subsec:estimativas_u_m_s<2} and \ref{subsec:estimativas_u_m_s_geq_2} into account, we test \eqref{v-u} by $v^n - \Delta v^n$ and conclude that
          \begin{equation} \vspace{-1pt}  
            v_m^k \mbox{ is bounded in } L^{\infty}(0, \infty; H^1(\Omega)),
            \label{limitacao_v_m^kr_1}
          \vspace{-1pt} \end{equation} 
          \begin{equation} \vspace{-1pt}  
            \Delta v_m^k \mbox{ is bounded in } L^2(0, \infty; L^2(\Omega)),
            \label{limitacao_Delta_v_m^kr}
          \vspace{-1pt} \end{equation} 
          \begin{equation} \vspace{-1pt}  
            \norma{\tilde{v}_m^k - v_m^k}{L^2(0,\infty;H^1(\Omega))} \leq C\, k.
            \label{limitacao_diferenca_v_m^k_v_m^kr}
          \vspace{-1pt} \end{equation} 
          Adapting the ideas used in Subsections \ref{subsec:estimativas_u_m_s<2} and \ref{subsec:estimativas_u_m_s_geq_2} to show that $z_m^k - \alpha$ is bounded in $L^2(0,\infty;H^2(\Omega))$ and that $\partial_t \tilde{z}_m^k$ is bounded in $L^2(0,\infty;L^{3/2}(\Omega))$, we prove that
          \begin{equation} \vspace{-1pt}  
            v_m^k \mbox{ is bounded in } L^2(0, \infty; H^2(\Omega)),
            \label{limitacao_v_m^kr_2}
          \vspace{-1pt} \end{equation}
          \begin{equation} \vspace{-1pt}  
            \partial_t \tilde{v}_m^k \mbox{ is bounded in } L^2(0, \infty; L^{3/2}(\Omega)).
            \label{limitacao_dt_v_m^k}
          \vspace{-1pt} \end{equation} 
          Then, using \eqref{limitacao_v_m^kr_1}, \eqref{limitacao_Delta_v_m^kr}, \eqref{limitacao_diferenca_v_m^k_v_m^kr}, \eqref{limitacao_v_m^kr_2}, \eqref{limitacao_dt_v_m^k} and Lemma \ref{lema_Simon}, we prove that there exists a function $v$ such that, up to a subsequence, we have
          \begin{equation} \vspace{-1pt}  
          \begin{array}{c}
            v_m^k \rightarrow v \mbox{ weakly* in } L^{\infty}(0,\infty;L^{\infty}(\Omega)) \cap L^{\infty}(0,\infty;H^1(\Omega)), \\
            v_m^k \rightarrow v \mbox{ weakly in } L^2(0,\infty;H^2(\Omega)), \\
            v_m^k \rightarrow v \mbox{ strongly in } L^2(0,T;H^1(\Omega)) \cap L^p(0,T;L^p(\Omega)), p \in [1, \infty), \mbox{ for each } T > 0, \\
            \mbox{and } \partial_t \tilde{v}_m^k \rightarrow \partial_t v \mbox{ weakly in } L^2(0,\infty;L^{3/2}(\Omega)).
          \end{array}
          \label{convergencia_v_m_s_intermediario}
        \vspace{-1pt} \end{equation} 
        
        Now, using the convergence obtained for $T^m(u_m^k)$ in Subsections \ref{subsec:estimativas_u_m_s<2} and \ref{subsec:estimativas_u_m_s_geq_2} and \eqref{convergencia_v_m_s_intermediario}, we conclude that the limit function $v$ is the unique solution of
        \begin{equation*} \vspace{-1pt}
          \partial_t v - \Delta v + u^s v = 0, \ \partial_\eta v \Big |_{\Gamma} = 0,
        \vspace{-1pt} \end{equation*}  
        where $u$ is the function of the pair $(u,z)$, fixed at the beginning of Subsection \ref{subsec:convergencia_de_v_m^kr}. Therefore, thanks to the uniqueness of the limit function $v$, we conclude that the whole sequence $v_m^k$ converges towards $v$ as $(m,k) \to (\infty,0)$. In addition, combining Lemma \ref{lema_equivalencia_problemas} and the uniqueness of the function $v$, given $u$, we deduce that $v = z^2 - \alpha^2$.
        
        Thus we conclude that $\{ (u_m^k, v_m^k) \}$ converges to $(u,v)$, a solution of \eqref{problema_P}, as $(m,k) \to (\infty, 0)$ in both cases, with $v^n$ given by \eqref{v-z} or \eqref{v-u}, finishing the proof of Theorem \ref{teo_principal}.


\section{Conclusion}

  In this work, we designed a conservative, energy stable, positivity preserving and convergent time discrete scheme for the chemotaxis-consumption model \eqref{problema_P}. Using the change of variables $z = \sqrt{v + \alpha^2}$ and a upper truncation of $u$ in the nonlinear chemotaxis and consumption terms, we proposed a Backward Euler scheme for the $(u,z)$-problem and two different ways of retrieving an approximation for $v$. We have proved the existence of solution to the time discrete scheme, uniform in time \emph{a priori} estimates and convergence of the scheme towards a weak solution $(u,v)$ of the chemotaxis-consumption model \eqref{problema_P}.

  We remark that, in the present work, we dealt with some issues that aroused from the time discretization of \eqref{problema_P}. Indeed, in order to obtain a time discrete scheme satisfying an energy law, independently of the time step size, it was crucial that the method was implicit and that the time discrete scheme was proposed in terms of the variable $z = \sqrt{v + \alpha^2}$ instead of the variable $v$. Moreover, Lemma \ref{lema_f''} was decisive to prove convergence in the case $s \in [1,2)$. We hope that this work will provide insights into how to propose a fully discrete scheme to approximate \eqref{problema_P}, rather as the time discretization performed in \cite{guillen2020study} contributed to the design of fully discrete schemes to approximate chemorepulsion-production models in \cite{guillen2020study2,guillen2021chemorepulsion}.

  Nevertheless, we expect that throughout the design of a conservative, energy stable, positivity preserving and convergent fully discrete scheme for the chemotaxis-consumption model \eqref{problema_P}, other difficulties may appear. In fact, as we can observe in the present work, the proof of energy stability (and hence convergence) for the proposed time discrete scheme relies on the cancellation between two terms which come from the chemotaxis and the consumption effects, respectively. We believe that this cancellation should also be the key to the energy stability of fully discrete schemes for \eqref{problema_P}. Since the cancellation depends strongly on the expression of the aforementioned terms, the design of an energy stable fully discrete scheme may possibly be quite challenging, especially if we want to preserve positivity.



\bibliographystyle{siamplain}
\bibliography{references}

\begin{thebibliography}{10}

\bibitem{bellomo2015toward}
{\sc N.~Bellomo, A.~Bellouquid, Y.~Tao, and M.~Winkler}, {\em Toward a
  mathematical theory of {K}eller--{S}egel models of pattern formation in
  biological tissues}, Mathematical Models and Methods in Applied Sciences, 25
  (2015), pp.~1663--1763.

\bibitem{ViannaGuillen2023uniform}
{\sc A.~L. Corr\^ea Vianna~Filho and F.~Guillén-González}, {\em Uniform in
  time solutions for a chemotaxis with potential consumption model}, Nonlinear
  Analysis: Real World Applications, 70 (2023), p.~103795.

\bibitem{duarte2021numerical}
{\sc A.~Duarte{-}Ro{\-}drí{\-}guez, M.~A. Rodríguez-Bellido, D.~A.
  Rueda-Gómez, and {\'E}.~J. Villamizar-Roa}, {\em Numerical analysis for a
  chemotaxis-{N}avier--{S}tokes system}, ESAIM. Mathematical Modelling and
  Numerical Analysis, 55 (2021), p.~417.

\bibitem{Evans2010}
{\sc L.~C. Evans}, {\em Partial {D}ifferential {E}quations}, Graduate Studies
  in Mathematics, American Mathematical Society, Providence, 2~ed., 2010.

\bibitem{eyre1998unconditionally}
{\sc D.~J. Eyre}, {\em An unconditionally stable one-step scheme for gradient
  systems}, Unpublished article,  (1998).

\bibitem{GilbargTrudinger}
{\sc D.~Gilbarg and N.~S. Trudinger}, {\em Elliptic {P}artial {D}ifferential
  {E}quations of {S}econd {O}rder}, Classics in Mathematics, Springer-Verlag,
  Berlin, 2001.

\bibitem{Grisvard}
{\sc P.~Grisvard}, {\em Elliptic {P}roblems in {N}onsmooth {D}omains},
  Monographs and studies in Mathematics, Pitman Publishing, Boston, 1985.

\bibitem{guillen2020study}
{\sc F.~Guill{\'e}n-Gonz{\'a}lez, M.~Rodr{\'\i}guez-Bellido, and D.~A.
  Rueda-G{\'o}mez}, {\em Study of a chemo-repulsion model with quadratic
  production. {P}art {I}: analysis of the continuous problem and time-discrete
  numerical schemes}, Computers \& Mathematics with Applications, 80 (2020),
  pp.~692--713.

\bibitem{guillen2020study2}
{\sc F.~Guill{\'e}n-Gonz{\'a}lez, M.~Rodr{\'\i}guez-Bellido, and D.~A.
  Rueda-G{\'o}mez}, {\em Study of a chemo-repulsion model with quadratic
  production. {P}art {II}: analysis of an unconditionally energy-stable fully
  discrete scheme}, Computers \& Mathematics with Applications, 80 (2020),
  pp.~636--652.

\bibitem{guillen2021chemorepulsion}
{\sc F.~Guill{\'e}n-Gonz{\'a}lez, M.~Rodr{\'\i}guez-Bellido, and D.~A.
  Rueda-G{\'o}mez}, {\em A chemorepulsion model with superlinear production:
  analysis of the continuous problem and two approximately positive and
  energy-stable schemes}, Advances in Computational Mathematics, 47 (2021),
  p.~87.

\bibitem{guillen2022comparison}
{\sc F.~Guill{\'e}n-Gonz{\'a}lez, M.~A. Rodr{\'\i}guez-Bellido, and D.~A.
  Rueda-G{\'o}mez}, {\em Comparison of two {F}inite {E}lement schemes for a
  chemo-repulsion system with quadratic production}, Applied Numerical
  Mathematics, 173 (2022), pp.~193--210.

\bibitem{guillen2019unconditionally}
{\sc F.~Guillén-González, M.~A. Rodríguez-Bellido, and D.~A. Rueda-Gómez},
  {\em Unconditionally energy stable fully discrete schemes for a
  chemo-repulsion model}, Mathematics of Computation, 88 (2019),
  pp.~2069--2099.

\bibitem{FranciscoGiordano}
{\sc F.~Guillén-González and G.~Tierra}, {\em Finite element numerical
  schemes for a chemo-attraction and consumption model}, arXiv preprint
  arXiv:2112.03431,  (2021).

\bibitem{jiang2015global}
{\sc J.~Jiang, H.~Wu, and S.~Zheng}, {\em Global existence and asymptotic
  behavior of solutions to a chemotaxis--fluid system on general bounded
  domains}, Asymptotic Analysis, 92 (2015), pp.~249--258.

\bibitem{keller1970initiation}
{\sc E.~F. Keller and L.~A. Segel}, {\em Initiation of slime mold aggregation
  viewed as an instability}, Journal of theoretical biology, 26 (1970),
  pp.~399--415.

\bibitem{keller1971model}
{\sc E.~F. Keller and L.~A. Segel}, {\em Model for chemotaxis}, Journal of
  theoretical biology, 30 (1971), pp.~225--234.

\bibitem{patlak1953random}
{\sc C.~S. Patlak}, {\em Random walk with persistence and external bias}, The
  bulletin of mathematical biophysics, 15 (1953), pp.~311--338.

\bibitem{Simon1986compact}
{\sc J.~Simon}, {\em Compact sets in the space ${L}^p (0, {T}; {B})$}, Annali
  di Matematica pura ed applicata, 146 (1986), pp.~65--96.

\bibitem{tao2011boundedness}
{\sc Y.~Tao}, {\em Boundedness in a chemotaxis model with oxygen consumption by
  bacteria}, Journal of mathematical analysis and applications, 381 (2011),
  pp.~521--529.

\bibitem{tao2012eventual}
{\sc Y.~Tao and M.~Winkler}, {\em Eventual smoothness and stabilization of
  large-data solutions in a three-dimensional chemotaxis system with
  consumption of chemoattractant}, Journal of Differential Equations, 252
  (2012), pp.~2520--2543.

\bibitem{winkler2012global}
{\sc M.~Winkler}, {\em Global large-data solutions in a
  {C}hemotaxis-({N}avier--) {S}tokes system modeling cellular swimming in fluid
  drops}, Communications in Partial Differential Equations, 37 (2012),
  pp.~319--351.

\end{thebibliography}
\end{document}